\documentclass[aop,preprint]{imsart}

\RequirePackage[OT1]{fontenc}
\RequirePackage{amsthm,amsmath,natbib}
\RequirePackage[colorlinks]{hyperref}
\RequirePackage{hypernat}

\numberwithin{equation}{section}
\usepackage{mathrsfs}
\usepackage{verbatim}
\usepackage{amssymb}
\usepackage{comment}
\usepackage{appendix}
\usepackage{graphicx}
\usepackage{subfigure}
\usepackage{enumerate}
\usepackage{comment}
\usepackage{pdfsync}
\usepackage{microtype}
\usepackage{subfigure}
\usepackage{color} 
\usepackage[table]{xcolor}
\usepackage{todonotes}
\usepackage{bbm, dsfont}
\usepackage{natbib}
\usepackage[left=3cm, right=3cm]{geometry}
\input xy
\xyoption{all}

\startlocaldefs

\newtheorem{theorem}{Theorem}[section]
\newtheorem{lemma}[theorem]{Lemma}
\newtheorem{proposition}[theorem]{Proposition}
\newtheorem{corollary}[theorem]{Corollary}

\newtheorem{definition}[theorem]{Definition}

\newtheorem{remark}[theorem]{Remark}

\newcommand{\E}{\mathbb{E}}

\newcommand*{\be}{\begin{equation}}
\newcommand*{\ee}{\end{equation}}
\newcommand*{\ba}{\begin{aligned}}
\newcommand*{\ea}{\end{aligned}}
\newcommand*{\barr}{\begin{array}{c}}
\newcommand*{\earr}{\end{array}}

\newcommand{\CMnd}{\mathrm{CM}_n(\boldsymbol{d})}
\newcommand{\indicator}{\mathbbm{1}}
\newcommand{\indic}[1]{\indicator_{\{#1\}}}
\newcommand{\indicwo}[1]{\indicator_{#1}}
\newcommand{\eqn}[1]{\begin{equation}#1\end{equation}}
\newcommand{\eqan}[1]{\begin{align}#1\end{align}}
\newcommand{\nn}{\nonumber}
\newcommand{\sss}{\scriptscriptstyle}
\newcommand{\e}{\mathrm{e}}

\newcommand{\Totan}{T^{\sss(1,2)}_{a_n}}
\newcommand{\Toan}{T^{\sss(1)}_{a_n/2}}
\newcommand{\Ttan}{T^{\sss(2)}_{a_n/2}}
\newcommand{\vep}{\varepsilon}
\newcommand{\SWG}{{\sf SWG}}

\newcommand{\whp}{whp }
\newcommand{\Nlos}{N_{\sss \mathcal{L}}}
\newcommand{\NLna}{N_{\sss \mathcal{L}_n}^*}
\newcommand{\NLnaa}{N_{\sss \mathcal{L}_n}^{**}}
\newcommand{\Nlaa}{N_{\sss \mathcal{L}}^{**}}

\newcommand{\prob}{\mathbb{P}}
\newcommand{\expec}{\mathbb{E}}

\newcommand{\convd}{\stackrel{d}{\longrightarrow}}
\newcommand{\convp}{\stackrel{\prob}{\longrightarrow}}

\definecolor{darkgreen}{rgb}{0,.4,0}
\definecolor{darkagenta}{rgb}{.5,0,.5}
\definecolor{darkred}{rgb}{1,0,0}
\definecolor{darkblue}{rgb}{0,0,.4}
\definecolor{black}{rgb}{0,0,0}

\newcommand{\Vlos}{V^{\sss (\mathcal{L})}}
\newcommand{\Vwin}{V^{\sss (\mathcal{W})}}

\newcommand{\tildeVwin}{\widetilde{V}^{\sss (\mathcal{W})}}
\newcommand{\tildeV}{\widetilde{V}}
\newcommand{\Xwin}{X^{\sss (\mathcal{W})}}

\newcommand{\Ltwin}{L_{\sss \mathcal{W}_n}^{\sss(t)}}
\newcommand{\barLtwin}{\bar{L}_{\sss \mathcal{W}_n}^{\sss(t)}}

\newcommand{\barLt}{\bar{L}_n^{\sss(t)}}
\newcommand{\Ntkwin}{N_{\sss \mathcal{W}_n}^{\sss(k,t)}}
\newcommand{\Ntkpluswin}{N_{\sss \mathcal{W}_n}^{\sss(k+1,t)}}

\newcommand{\barNtkwin}{\bar{N}_{\sss \mathcal{W}_n}^{\sss(k,t)}}
\newcommand{\barNtkpluswin}{\bar{N}_{\sss \mathcal{W}_n}^{\sss(k+1,t)}}
\newcommand{\barNtk}{\bar{N}^{\sss(k,t)}_n}

\newcommand{\field}{\mathscr{F}}
\newcommand{\op}{o_{\sss \prob}}
\newcommand{\Op}{O_{\sss \prob}}
\newcommand{\Thetap}{\Theta_{\sss \prob}}
\newcommand{\Ver}{U}

\newcommand{\Tlos}{T^{\sss({\mathcal{L}})}}

\newcommand{\BPn}{{\sf BP}_n}
\newcommand{\YBP}{Y^{\sss\rm({\sf BP})}}
\newcommand{\YSWG}{Y^{\sss\rm({\sf SWG})}}
\newcommand{\Alive}{{\sf A}}
\newcommand{\GF}{F_{\sss X}}

\newcommand{\ra}{\rightarrow}

\DeclareSymbolFont{extraup}{U}{zavm}{m}{n}
\DeclareMathSymbol{\varheart}{\mathalpha}{extraup}{86}
\DeclareMathSymbol{\vardiamond}{\mathalpha}{extraup}{87}
\newcommand{\ensymboldefinition}{\blacktriangleleft}

\begin{document}
\begin{frontmatter}
\title{Universal `winner-takes-it-all' phenomenon in~scale-free~random~graphs}
\runtitle{Universal winner-takes-it-all}

\begin{aug}

\runauthor{R. van der Hofstad}

\author{Remco van der Hofstad}
\address{Eindhoven University of Technology,\protect\\
\tt{r.w.v.d.hofstad@TUE.nl}}


\affiliation{Eindhoven University of Technology}

\end{aug}
\date{\today}

\begin{abstract} 
We study competition on scale-free random graphs, where the degree distribution satisfies an asymptotic power-law with infinite variance. Our competition process is such that the two types attempt at occupying vertices incident to the presently occupied sets, and the passage times are independent and identically distributed, possibly with different distributions for the two types. Once vertices are occupied by a type, they remain on being so. 

We focus on the explosive setting, where our main result shows that the winning type occupies all but a {\em finite} number of vertices. This universal `winner-takes-it-all' phenomenon significantly generalises previous work with Deijfen for exponential edge-weights, and considerably simplifies its proof.  
\end{abstract}



\end{frontmatter}

\maketitle


\section{Introduction}
\label{sec-intro}

We study {\em competition problems} on scale-free random graphs, where edges are being traversed by two different types and vertices remain on having the type they received upon their first occupation. We assume that the edge-traversal times are independent and identically distributed as in first-passage percolation, possibly with different distributions for the two types. Thus, the competition processes each follow the first-passage percolation spreading rules, but the exclusion of vertices that are already occupied by the other type(s) slows the spreading down, effectively creating {\em competition} between the types. 

We consider the configuration model, one of the most popular random graph models having flexible degree distributions, in the scale-free setting where the degrees have {\em infinite variance}. Prior work with Deijfen \cite{DeiHof16} investigates this in the (explosive) setting where the edge-traversal times of the competition have exponential distributions with possibly different parameters. Baroni, van der Hofstad, and Komj\'athy \cite{BarHofKom15} and Komj\'athy \cite{Kom15} study this problem in the graph-distance setting where the edge-traversal times are deterministic, with either different \cite{BarHofKom15} or equal \cite{Kom15} deterministic traversal times. Our work significantly extends that with Deijfen \cite{DeiHof16}.

The main issue in competition problems is to determine which of the types occupies the majority of the graph, or whether both can occupy a substantial portion of it. The latter case is called {\em co-existence}, and so far has been proved to occur with positive probability only for deterministic and equal edge-weights in the scale-free random graph setting \cite{Kom15}. In  \cite{DeiHof16}, it is proved that for exponential edge-weights, the losing type occupies a number of vertices that converges in distribution. It is this result that we aim to extend from the Markovian setting with exponential edge-weights to general edge-weights that are {\em explosive}, in that first-passage percolation with such edge-weights reaches an unbounded number of vertices in a tight amount of time. This explosive setting only occurs for configuration models with infinite-variance degrees. The setting of competition on configuration models with finite-variance degrees was studied by Ahlberg, Deijfen and Janson \cite{AhlDeiJan17} in the Markovian setting, who show that coexistence occurs precisely when the edge-traversal times have the same means. With Deijfen and Sfragara, we investigate the setting of configuration models with {\em infinite-mean degrees} in \cite{DeiHofSfr22}, where the losing type is likely to only occupy its initial starting point. See \cite{HofShn23} for a related problem aimed at modelling the spread of misinformation.

The closely related problem of first-passage percolation on random graphs has attracted considerable attention, see \cite[Chapter 3]{Hofs23} for an overview. This line of research started with joint work with Bhamidi and Hooghiemstra \cite{BhaHofHoo12b}, where it was shown that the fluctuations of the weight of, and number of edges in, optimal paths in first-passage percolation on random graphs with finite-variance degrees behave remarkably universally, extending work in \cite{BhaHofHoo09b} and \cite{BhaHofHoo12} that applies to exponential edge-weights and finite-variance degrees. The case with infinite-variance degrees (investigated for exponential edge-weights in \cite{BhaHofHoo09b} as well) is rather different, since there are several universality classes depending on whether the age-dependent (continuous-time) branching process that approximates the flow through the random graph, is explosive or non-explosive, which we refer to as being {\em conservative}. In the explosive setting, the weight-distance between two uniform vertices converges in distribution to the sum of two independent explosion times. In the conservative setting, on the other hand, the weight-distance between two uniform vertices converges to infinity in probability \cite{BarHofKom17}. See also \cite{BarHofKom19} where a particular conservative setting is explored in more detail.

\paragraph{\bf Main result and innovation} The first main result in this paper shows that in scale-free configuration models where at least one of the edge-weight distributions is explosive, one of the two types occupies the majority of the graph. We refer to this type as the {\em winning} type, and the other as the {\em losing} type. In particular, co-existence occurs with vanishing probability. The second main result shows that the losing type occupies a {\em finite} number of vertices, and we identify its limiting distribution. This was done first for exponential edge-weights in \cite{DeiHof16}. Our proof is simpler than that in \cite{DeiHof16}, and relies on a key first-passage percolation result in that we identify the {\em epidemic curve}, which is interesting in its own right. When viewing first-passage percolation as the spread of a disease, the epidemic curve describes the proportion of vertices that have contracted the disease as a function of time. 

For the explosive settings that we investigate, we can think of the time of explosion as the time at which one of the highest-degree vertices, or {\em hubs}, is found. Our epidemic curve result shows that at any time $t>0$ after explosion, a positive proportion of the vertices is occupied by the first-passage percolation. In particular, vertices of high degree are more likely to have been occupied. When lifting these results to the competition setting, this means that virtually immediately after the time that the winning type occupies a hub, it has found all high-degree vertices, thus dramatically limiting the possible occupation of vertices by the losing type. In fact, we prove that immediately after explosion of the winning type, the exploration of the losing type in the part of the graph that has not yet been occupied by the winning type becomes {\em conservative}.  In turn, this implies that the losing type only occupies a tight number of vertices that converges in distribution. In this last step, we generalise a technical bound obtained in \cite{DeiHof16} to general explosive edge-weights.

\paragraph{\bf Notation and organization}
We abbreviate {\em left}- and {\em right-hand side} by lhs and rhs, respectively. We abbreviate {\em uniformly at random} by uar, and {\em with respect to} by wrt.
Given two $\mathbb{R}$-valued random variables $X$ and $Y$, we say that $X\stackrel{d}{=}Y$ if  $\prob(X\leq x)=\prob(Y\leq x)$ for all $x\in \mathbb{R}$. A sequence of random variables $(X_n)_{n\geq 1}$ {\em converges in probability} to a random variable $X$, denoted as $X_{n}\stackrel{\prob}{\rightarrow}X$, if, for all $\varepsilon>0$, $\prob(|X_{n}-X|>\varepsilon)\rightarrow 0$, and it {\em converges in distribution} to $X$, denoted as $X_{n}\stackrel{d}{\rightarrow}X,$ if $\lim_{n\rightarrow\infty}\prob(X_{n}\leq x)=\prob(X\leq x)$ for all $x\in \mathbb{R}$ for which $F_{\sss{X}}(x)=\prob(X\leq x)$ is continuous. A sequence of random variables $(X_{n})_{n\geq 1}$ is said to be {\em tight}, which we denote as $X_n=\Op(1)$, if for all $\varepsilon>0$, there exists $r>0$ such that $\sup_{n\geq 1}\prob(|X_{n}|>r)<\varepsilon$. We write $X_n=\op(1)$ when $X_n\convp 0$. 
A sequence of events $(\mathcal{E}_{n})_{n\geq 1}$ is said to hold {\em with high probability} (whp) if $\lim_{n\rightarrow\infty}\prob(\mathcal{E}_{n})=1$. 

We denote the set $\{1,2,\dots, n\}$ by $[n]$. $\CMnd$ denotes the configuration model on $n$ vertices with independent and identically distributed (iid) degrees $\boldsymbol{d}=(d_v)_{v\in[n]}$, as introduced formally in the next section.
 
This section is organised as follows. In Section \ref{sec-res}, we introduce the model and state our results.  We close this section in Section \ref{sec-overview} by giving an overview of the proof, and in Section \ref{sec-disc} by relating it to the literature, as well as giving some open problems.

\subsection{Model and results}
\label{sec-res}
Our setting is the configuration model $\CMnd$ originally defined by Bollob\'as \cite{Boll80b} for random regular graphs, and studied in great detail by Molloy and Reed \cite{MolRee95, MolRee98}.  $\CMnd$ has $n$ vertices, where $\boldsymbol{d}=(d_1,\dots, d_n)$, and $d_v$ is the degree of vertex $v\in[n]$. Further, let $L_{n}=\sum_{v\in[n]}d_v$ denote the total degree. The configuration model random graph is obtained as follows: We assign a number $d_v$ of half-edges to vertex $v$, and we pair these half-edges uniformly at random. When two half-edges are paired they form an edge and are removed from the set of unpaired half-edges. We continue the procedure until there are no more unpaired half-edges available.
If $L_{n}$ is odd, then we add a half-edge to vertex $n$; this extra half-edge makes no difference to our results and we will refrain from discussing it further. We denote the edge set of $\CMnd$ by $E(\CMnd)$. We consider the configuration model $\CMnd$ on $n$ vertices with iid degrees. We assume that the degrees have a power-law distribution, i.e., for some constants $c_1,C_1$ and power-law exponent $\tau$ with $\tau\in(2,3)$,
	\be
	\label{F-conditions}
	\frac{c_1}{x^{\tau-1}}\le 1- F_{\sss D}(x)= \prob(D>x) \le \frac{C_1}{x^{\tau-1}}.
	\ee
In this case, the degrees have finite mean, but infinite variance. We further assume that
	\be
	\label{gc}
	\prob(D\geq 2)=1,
	\ee 
which is equivalent to the largest connected component consisting of most of the vertices (see e.g., \cite{JanLuc09, BolRio15, Hofs21,Hofs24}, or \cite{FedHof17, Fede20} for the sharpest conditions). This describes the random graph on which our competition process runs.

We next describe the edge-traversal times, often called {\em edge-weights}, that characterise the distribution of the spread of the competition along the edges.  We associate to each edge two distinct and independent families of iid  edge-weights denoted by $(X^{\sss(1)}_e)_{e\in E(\CMnd)}$ and $(X^{\sss(2)}_e)_{e\in E(\CMnd)}$. These sequences are each iid and independent from each other, but $(X^{\sss(1)}_e)_{e\in E(\CMnd)}$ and $(X^{\sss(2)}_e)_{e\in E(\CMnd)}$ could have different distributions. We do assume that at least one of the edge-weight distributions $X^{\sss(1)}$ and $X^{\sss(2)}$ is such that the age-dependent branching process that locally describes $\CMnd$ with its edge-weighted geometry is explosive. For this, let $D^\star-1$ be defined as the (\emph{size-biased version} of $D$)$-1$, that is,
	\begin{align}
	\label{forwarddegree}
	\prob(D^\star-1=k):=\frac{k+1}{\E[D]}\prob(D=k+1).
	\end{align}
We let $F^\star$ denote the distribution function of $D^\star-1$. Then, our age-dependent branching process has offspring distribution $D^\star-1$. We call an age-dependent branching process with edge-weight (or life-time) distribution $X$ and offspring distribution $D^\star-1$ a {\em $(D^\star-1,X)$ age-dependent branching process.} When the root has offspring distribution $D$, while every other individual has offspring distribution $D^\star-1$, we call this process the {\em unimodular $(D,X)$ age-dependent branching process}. Conditions for explosion of this age-dependent branching processes are discussed in \cite{BarHofKom17}, as well as in the references therein. See also \cite{Komj16} for a detailed discussion of explosion times of continuous-time branching processes for epidemics, including contagious periods and infection times.

We consider the spread of two {\em independent} infections from the vertices $\Ver_1$ and $\Ver_2$, where $X^{\sss(i)}_e$, with $i\in\{1,2\},$ is the passage cost of edge $e$ and $\Ver_1, \Ver_2$ are chosen uar and independently from $[n]$. Note that, since the degrees in $\CMnd$ are iid, the vertices in the configuration model are {\em exchangeable}, so that we might as well have chosen vertices 1 and 2 as the starting locations for our competitions. 

Considering the infection as a spread of colour in the graph, at time $0$, the vertex $\Ver_1$ is infected by red, and vertex $\Ver_2$ is infected by blue.
Each infection spreads through the edges, the time that the red and blue take to go through the edge $e$ is given by $X^{\sss(1)}_e$, and $X^{\sss(2)}_e,$ respectively. We assume that once a vertex is infected by a colour, it remains of that colour forever.

Our main result concerns the distribution of the winning type, as well as the number of vertices occupied by the losing type:

\begin{theorem}[Universal ``winner-takes-it-all'' phenomenon]
\label{thm-main} 
Consider $\CMnd$ with iid degrees satisfying \eqref{F-conditions}--\eqref{gc}.
Assume that at least one of the age-dependent branching processes with offspring distribution $D^\star-1$ and edge-weights $X^{\sss(1)}_e$, and $X^{\sss(2)}_e,$ respectively, is explosive. Denote the respective explosion times of the unimodular $(D,X^{\sss(i)})$ age-dependent branching processes by $V^{\sss(1)}$ and $V^{\sss(2)}$, respectively (where at most one of these can be infinite a.s.). Then
\begin{itemize}
\item[\rm{(a)}] the fraction $\bar{N}_1(n)$ of type 1 infected vertices converges in distribution to the indicator variable $\indic{V^{\sss(1)}<V^{\sss(2)}}$ as $n\to\infty$;
\item[\rm{(b)}] the total number $\Nlos(n)$ of vertices occupied by the losing type converges in distribution to a proper random variable $\Nlos$.
\end{itemize}
\end{theorem}

The main message in Theorem \ref{thm-main} is that the losing type only occupies a {\em finite} number of vertices, and that this number converges in distribution to a limiting proper random variables: a universal ``winner-takes-it-all'' phenomenon (recall \cite{DeiHof16}). When the unimodular $(D,X^{\sss(1)})$ and $(D,X^{\sss(2)})$ branching processes are {\em both} explosive, it is not hard to see that $\prob(V^{\sss(1)}<V^{\sss(2)})\in(0,1)$, so each of the two types can win the majority of the vertices with positive probability. Thus, a type that is explosive has a positive probability of being the losing type, but a type that is conservative always loses from the explosive type. This is quite different from the Markovian setting of configuration models with finite-variance degrees, as studied by Ahlberg, Deijfen and Janson \cite{AhlDeiJan17} for configuration models with general finite-variance degrees, and by Antunovic et al.\ \cite{AntDekMosPer11} for random regular graphs, and can be understood by noting that such age-dependent branching processes are never explosive.

A classical example to which our results apply is when $X^{\sss(2)}$ has the same law as $\lambda X^{\sss(1)}$ for some $\lambda$ (as for example in \cite{DeiHof16}), in which case $V^{\sss(2)}$ has the same distribution as $\lambda V^{\sss(1)}$. Thus, the unimodular $(D,X^{\sss(1)})$ and $(D,X^{\sss(2)})$ branching processes just differ in the {\em speed} at which the types occupy the vertices. Yet, the slowest type can still occupy the majority of the vertices. 

\begin{remark}[Other sources for the competition]
\label{rem-starting-points}
{\rm We prove Theorem \ref{thm-main} only in the case where the starting points are each chosen uar from $[n]$, independently of each other. However, the proof reveals that a similar statement also holds in the case where the starting points are chosen as the ends of a uniform {\em edge}. In this case, the explosion times $V^{\sss(1)}$ and $V^{\sss(2)}$ of the unimodular $(D,X^{\sss(1)})$ and $(D,X^{\sss(2)})$ branching processes are replaced by the explosion times $\widetilde V^{\sss(1)}$ and $\widetilde V^{\sss(2)}$ of the $(D^\star-1,X^{\sss(1)})$ and $(D^\star-1,X^{\sss(2)})$ age-dependent branching processes, where also the root has offspring $D^\star-1$.}\hfill$\ensymboldefinition$
\end{remark}

\subsection{Overview of the proof}
\label{sec-overview}
The keys steps in the proof are as follows:

\paragraph{\bf Stage up to first explosion} 
We start by growing the neighborhoods of the vertices $\Ver_1,\Ver_2$ in the first-passage percolation setting. Fix $\rho>0$ appropriately. By \cite[Proposition 4.7]{BhaHofHoo09b}, the degrees of the vertices found up to time $n^{\rho}$ can be coupled to  a collection of $a_n=n^\rho$ iid random variables having distribution $D^\star-1$, except for the degrees of $\Ver_1,\Ver_2$, which are independent copies of the asymptotic degree distribution $D$. Denote the smallest-weight graph thus obtained when exploring up to time $t$ by $\SWG^{\sss(1,2)}(t)$  (so that $\SWG^{\sss(1,2)}(0)=\{\Ver_1,\Ver_2\}$).  Let $T_m^{\sss(1,2)}$ denote the time at which the $(m+2)$nd vertex was found (so that $T_0^{\sss(1,2)}=0$). 
Then, again by \cite[Proposition 4.7]{BhaHofHoo09b}, for $\rho>0$ sufficiently small, whp $\SWG^{\sss(1,2)}(T_{a_n}^{\sss(1,2)})$ consists of two disjoint sets of vertices, one containing $\Ver_1$, the other containing $\Ver_2$. We introduce such smallest-weight graphs in more detail in Section \ref{sec-SWGs}, following \cite{BhaHofHoo09b, BhaHofKom14}, as well as \cite{BhaHofHoo12b}.

Since the edge-weights  $(X^{\sss(1)}_e)_{e\in E(\CMnd)}$ and $(X^{\sss(2)}_e)_{e\in E(\CMnd)}$ are independent, this means that $(\SWG^{\sss(1,2)}(T_m^{\sss(1,2)}))_{m\leq a_n}$ is effectively coupled to the disjoint union of two independent unimodular $(D, X^{\sss(1)})$ and $(D,X^{\sss(2)})$ branching processes. Let $T_m^{\sss(i)}$ denote the time at which the $(m+1)$st vertex was found in the first-passage exploration from vertex $\Ver_i$, so that $(T_m^{\sss(1,2)})_{m\leq a_n}$ is the ordered sequence of $(T_m^{\sss(1)})_{m\leq a_n}$ and $(T_m^{\sss(2)})_{m\leq a_n}$ truncated to length $a_n$. Note that $T_{a_n}^{\sss(1,2)}\geq T_{a_n/2}^{\sss(1)}\wedge T_{a_n/2}^{\sss(2)}$, where $x\wedge y$ denotes the minimum of $x,y\in {\mathbb R}$. In our proof, it turns out to be more convenient to work at time $T_{a_n/2}^{\sss(1)}\wedge T_{a_n/2}^{\sss(2)}$ instead of $T_{a_n}^{\sss(1,2)}$, since this allows for a more convenient use of the asymptotic independence of $(\SWG^{\sss(1)}(T_m^{\sss(1)}))_{m\geq 0}$ and $(\SWG^{\sss(2)}(T_m^{\sss(2)}))_{m\geq 0}$. Since, for $a_n\rightarrow \infty$ sufficiently slowly, $(T_{a_n/2}^{\sss(1)}, T_{a_n/2}^{\sss(2)})\convd (V^{\sss(1)},V^{\sss(2)})$, and $\prob(V^{\sss(1)}\neq V^{\sss(2)})=1$, it means that one of the trees in $\SWG^{\sss(1,2)}(T_{a_n/2}^{\sss(1)}\wedge T_{a_n/2}^{\sss(2)})$ has size $a_n/2-\Op(1)$, while the other has size $\Op(1)$, the larger tree corresponding to the smaller $T_{a_n/2}^{\sss(i)}$. 

Without loss of generality, we may assume that $T_{a_n/2}^{\sss(1)}<T_{a_n/2}^{\sss(2)}$. Then, type 1 is quite ahead of type 2 in the competition process. We will show that whp vertex 1 then wins the majority of the vertices. The above coupling to continuous-time branching processes forms the starting point of the analysis, in which the competition did not affect the spread of types 1 and 2 too much. This changes dramatically immediately afterwards, as we explain next.

\paragraph{\bf Slightly beyond the first explosion} 
Let $M_n^{\sss (i)}$ denote the size of the tree containing $\Ver_i$ at time $T_{a_n/2}^{\sss(1)}\wedge T_{a_n/2}^{\sss(2)}$, so that $M_n^{\sss (2)}\convd M^{\sss (2)}$ for some finite integer-valued random variable $M^{\sss (2)}$ on the event that $T_{a_n/2}^{\sss(1)}<T_{a_n/2}^{\sss(2)}$. Since type 2 explodes later than type 1, and the inter-arrival times $T_{m+1}^{\sss(2)}-T_m^{\sss(2)}$ are continuous random variables, $T_{M_n^{\sss (2)}+1}^{\sss(2)}-T_{M_n^{\sss (2)}}^{\sss(2)}>\vep$ whp for small $\vep$. Thus, for some time after time $T_{a_n/2}^{\sss(1)}\wedge T_{a_n/2}^{\sss(2)}$, type 2 does not even attempt to occupy another vertex. Type 1, on the other hand, is hardly hindered by the competition due to type 2 (as there are a tight number of vertices occupied by type 2), and thus proceeds as a normal first-passage percolation flow after time $T_{a_n/2}^{\sss(1)}\wedge T_{a_n/2}^{\sss(2)}$. At time $T_{a_n/2}^{\sss(1)}\wedge T_{a_n/2}^{\sss(2)}+\vep,$ this flow has occupied some $\Thetap(n) f(\vep)$ vertices, where $f(\vep)>0$ for every $\vep>0$. Moreover, vertices with high degree are more likely to have been occupied by type 1 when $T_{a_n/2}^{\sss(1)}<T_{a_n/2}^{\sss(2)}$ (see Proposition \ref{prop-bd-sec-mom}). Further, the two SWGs at time $T_{a_n/2}^{\sss(1)}\wedge T_{a_n/2}^{\sss(2)}+\vep$ are still disjoint whp (see Proposition \ref{prop-disj-trees}). 

\paragraph{\bf Exploration curve winning type beyond the first explosion} 
We conclude that immediately after the explosion of the winning type, the losing type can hardly make use of high-degree vertices (see Proposition \ref{prop-bd-sec-mom}), so that its spread becomes conservative (see Lemma \ref{lem-bd-growth-losing}). This leaves very little room for type 2 to occupy other vertices. As a result, the epidemic curve for the winning type is identical to that of first-passage percolation in the absence of competition, as formulated for first-passage percolation in Theorem \ref{thm-degree-time-FPP}, and for the competition process in Proposition \ref{prop-degree-time}. In particular, this shows that $N_{\sss \rm{los}}(n)=\op(n)$.

We prove that the above implies that $N_{\sss \rm{los}}(n)\convd N_{\sss \rm{los}}$ for some a.s.\ {\em finite} limiting variable $N_{\sss \rm{los}}$ by first showing that the proportion of vertices of degree $k$ found by the winning type at time $T_{a_n/2}^{\sss(1)}\wedge T_{a_n/2}^{\sss(2)}+t$ converges to some deterministic curve (see Proposition \ref{prop-degree-time}). It relies on the previous step that shows that the loosing type hardly influences the occupation process of the winning type. This now identifies the random environment in which the losing type is left to occupy vertices. It turns out that this environment is so restrictive, that the losing type only occupies a {\em finite} number of vertices, and its distribution can be characterised in terms of the exploration and the shape of the epidemic curve. This concludes our overview of the proof of our main result Theorem \ref{thm-main}.

\subsection{Discussion and open problems}
\label{sec-disc}

\paragraph{\bf Comparison with the proof in \cite{DeiHof16}} The proof strategy in this paper crucially relies on the epidemic curve result in Theorem \ref{thm-degree-time-FPP}, and its extension to competition processes in Proposition \ref{prop-degree-time}. As such, some of the technical computations and bounds in \cite{DeiHof16} for exponential edge-weights could be avoided. This both generalises as well as simplifies the proof.

\paragraph{\bf Multiple starting points and dependence} In \cite{DeiHof16}, also the setting where both types start from several independently chosen uniform locations was studied. Let $k_1$ and $k_2$ denote the number of starting points of types 1 and 2, respectively. An analysis was performed in the case where $k_1=k$ and $k_2=1$, as $k\rightarrow \infty$. This analysis basically boils down to the question how $\min_{j\leq k} V_j^{\sss(1)}$ behaves. It would be of interest to extend this analysis to our more general setting. Another interesting extension involves the setting where $(X^{\sss(1)}_e,X^{\sss(2)}_e)_e$ are iid, but $X^{\sss(1)}_e$ and $X^{\sss(2)}_e$ are not necessarily independent. We conjecture that the results in this paper are true in this setting as well.

\paragraph{\bf Conservative settings} It would be of interest to extend our analysis to the setting where {\em both} unimodular $(D,X^{\sss(1)})$ and $(D,X^{\sss(2)})$ branching processes are conservative. Let us for the moment stick to the setting where $X^{\sss(2)}$ has the same law as $\lambda X^{\sss(1)}$ for some $\lambda$. In this case, we believe that $T_{a_n}^{\sss(1)}$ and $\lambda T_{a_n}^{\sss(2)}$ grow at roughly the same {\em deterministic} speeds. When $\lambda>1$, this implies that type 1 occupies the majority of the vertices whp, and the fastest type wins. When $\lambda=1$, however, this suggests that the size of the smallest-weight graphs $M_n^{\sss (1)}$ and $M_n^{\sss (2)}$ at time $T_{a_n/2}^{\sss(1)}\wedge T_{a_n/2}^{\sss(2)}$ are {\em both} close to $a_n/2$. This raises the question whether co-existence can occur in this setting. The case where the passage times are deterministic has been addressed in  \cite{BarHofKom15, Kom15}, and shows that co-existence can indeed occur, the proof of which is remarkably subtle. 

\paragraph{\bf Finite-variance degrees} The setting of competition on configuration models with finite-variance degrees was carried out in \cite{AhlDeiJan17} for the Markovian case of exponential edge-weights (see also \cite{AntDekMosPer11}). It would be of great interest to extend such results to general edge-weights, where also the edge-weights of the two types might have a different distribution. The universality of first-passage percolation on configuration models, as proved in \cite{BhaHofHoo12b}, should be extremely helpful in this analysis. One may conjecture that co-existence can only occur when the Malthusian parameters of the corresponding CTBPs agree, but it is not unlikely that high-order corrections of the $n$-dependent Malthusian parameters play a role.

\paragraph{\bf Other random graphs} We believe that our results can be extended to several other random graphs with infinite variance degrees. Natural examples include preferential attachment models, where the first-passage percolation problem was investigated in \cite{JorKom20}, and (finite-graph versions of) scale-free percolation (see \cite{DeiHofHoo13}, and \cite{HofKom17b} for first-passage percolation results in this setting).

\paragraph{\bf Organization of the paper} In Section \ref{sec-FPP-prelim}, we recall previous results for first-passage percolation, as well as state and prove the epidemic curve result. In Section \ref{sec-most}, we show that the type that first grows large occupies all but a vanishing proportion of the graph. In Section \ref{sec-almost-all}, we extend this analysis by identifying the limit in distribution of the number of vertices occupied by the losing type.

\section{First-passage percolation preliminaries}
\label{sec-FPP-prelim}
In this section, we present some results concerning first-passage percolation on scale-free configuration models that are crucial in the proof of Theorem \ref{thm-main}. This section is organised as follows. In Section \ref{sec-SWGs}, we give more detailed information on first-passage percolation explorations, since these are crucial in the remainder of this paper. In Section \ref{sec-epidemic-curve}, we state a result that describes how first-passage percolation quickly sweeps through the graph right after explosion, see Theorem \ref{thm-degree-time-FPP} below. This result is interesting in its own right. In Section \ref{sec-prel}, we state further more technical preliminaries on first-passage percolation, which we use in Section \ref{sec-proof-epidemic-curve} to complete the proof of Theorem \ref{thm-degree-time-FPP}. In Section \ref{sec-FPP-prelim}, we only deal with first-passage percolation, and write the edge-weights as $(X_e)_{e\in E(\CMnd)}$, which are iid continuous random variables, whose unimodular $(D,X)$ branching process is explosive.

\subsection{Smallest-weight graphs for first-passage percolation explorations}
\label{sec-SWGs}
In this section, we discuss smallest-weight graphs, used for first-passage percolation explorations. We define our graph-exploration process in {\em continuous time}, and let $\SWG(t)$ denote the collection of vertices that can be reached within time $t$. We first describe this in general, and then specialise to the configuration model, for which we can perform the graph exploration simultaneously with the first-passage percolation exploration. Let $G=(V(G),E(G))$ be a general graph, and let $(X_e)_{e\in E(G)}$ denote the edge-weights. 

\paragraph{First-passage percolation graph exploration} Let $v\in V(G)$ be the source of our exploration, and set $\SWG(0)=\{v\}$. We call $N(v)$, the edges incident to $v$, {\em active}, and denote the active set of edges ${\sf A}(0)=N(v)$. We then let
	\eqn{
	T_1=\min\{X_{e}\colon v\in e\}
	}
be the first time at which the $\SWG(t)$ increases by one vertex. Let $v_1$ be such that $T_1=X_{\{v,v_1\}}$. Then, $\SWG(t)=\{v\}, {\sf A}(t)={\sf A}(0)$ for all $t\in [0,T_1)$, and $\SWG(T_1)=\{v,v_1\}$, while ${\sf A}(T_1)$ consists of all the edges $(N(v)\cup N(v_1))\setminus \{e_1\}$ incident to $v$ and $v_1$, except for the edge $e_1=\{v,v_1\}$. Finally, we let $W(v,v_1)=T_1$. We iterate this procedure, and let $\SWG(T_m)=\{v,v_1,\ldots, v_m\}$, and denote the set of active edges at time $T_m$ by ${\sf A}(T_m)$, which consists of all edges incident to the vertices $\{v,v_1,\ldots, v_m\}$, except for the edges $e_1, \ldots, e_m$, where $e_m=\{u_m,v_m\}$ is incident to $v_m$, while $u_m\in \SWG(T_{m-1})$. Finally, for every $u\in \SWG(T_m)$, $W_n(v,u)$ is the weighted graph distance or traversal time between $u$ and $v$. We then let
	\eqan{
	\label{T-m+1}
	T_{m+1}&=\min\{X_{\{u,v'\}}+W_n(v,u)\colon u\in \SWG(T_m), v\not\in \SWG(T_m), \{u,v\}\in E(G)\}\\
	&=\min\{X_{\{u,v'\}}+W_n(v,u)\colon \{u,v'\}\in {\sf A}(T_m), u\in \SWG(T_m)\},\nn
	}
we let $e_{m+1}=\{u_{m+1},v_{m+1}\}$ be such that $T_{m+1}=X_{\{u_{m+1},v_{m+1}\}}+W_n(v,u_{m+1})$, so that $e_{m+1}$ is the minimiser of \eqref{T-m+1}. Then, we set $W_n(v,v_{m+1})=T_{m+1}$, $\SWG(t)=\{v,v_1,\ldots, v_{m}\}$ and ${\sf A}(t)={\sf A}(T_m)$ for all $t\in [T_m,T_{m+1})$, while $\SWG(T_{m+1})=\{v,v_1,\ldots, v_{m+1}\}$ and let ${\sf A}(T_{m+1})$ be $({\sf A}(T_m)\cup N(v_{m+1}))\setminus \{e_{m+1}\}$ with all the edges between $N(v_{m+1})$ and ${\sf A}(T_m)$ removed (if present). This defines the first-passage percolation process in general. 
	 
\paragraph{First-passage percolation configuration model exploration}	 
We now specialise to the configuration model, for which we construct the graph {\em simultaneously} with the first-passage percolation exploration. Indeed, the configuration model can be constructed by successively pairing the half-edges. In our construction, it will be convenient to consider a set of active {\em half-edges} rather than the set of active edges. However, since in \eqref{T-m+1}, we only consider edges $\{u,v\}$ that are incident to $\SWG(T_m)$, i.e., we ignore minimal edges that give rise to cycles, we need to {\em screen} the half-edges depending on whether they give rise to cycles or not. For this, we let $v_{m+1}\in [n]$ be the minimiser in \eqref{T-m+1}, and consider the set of {\em possible} active half-edges to be all the half-edges incident to $v$ except for the half-edge that is part of the edge $e_{m+1}$ between $v_{m+1}$ and $\SWG(T_{m+1})$. We then, one by one, check whether these half-edges form cycles or self-loops. Those that do, together with the half-edges to which they are paired, are removed from the set of possible active half-edges, to obtain the set of active half-edges ${\sf A}(0)$. We call this the {\em cycle-removal procedure}. We give each of these half-edges an {\em independent} edge-weight. 

Note that rather than associating the edge-weights to the {\em edges}, we now associate them to {\em half-edges}. This has the risk that an edge has {\em two} edge-weights associated to it. It is for this reason that we need to remove edges that close cycles from the set of active edges. For a half-edge $y$, we let $v_y$ be the vertex incident to it. We then define, as in \eqref{T-m+1},
	\eqan{
	\label{T-m+1-CM}
	T_{m+1}&=\min\{X_{a}+W_n(v,v_a)\colon a\in {\sf A}(T_m)\},\nn
	}
where $a$ is a half-edge, and $v_a\in \SWG(T_m)$ is the vertex incident to half-edge $a$. After identifying the minimising half-edge, we pair it uniformly to a half-edge $y_{m+1}$ from the set of available and inactive half-edges. For all other half-edges incident to $v_{y_{m+1}}$, i.e., those unequal to $y_{m+1}$, we check whether they are paired to an active half-edge, or to another half-edge incident to $y_{m+1}$. Then, ${\sf A}(T_{m+1})$ consists of ${\sf A}(T_m)$ together with the half-edges incident to $v_{m+1}$, with the half-edge $y_{m+1}$ and the half-edges incident to $v_{y_{m+1}}$ that produce cycles and self-loops, as well as their pairs, removed. Due to this cycle-removal procedure, in fact all active half-edges in ${\sf A}(t)$ are connected to vertices {\em outside} of $\SWG(t)$, as in \eqref{T-m+1}.

\begin{remark}[Related first-passage percolation explorations]
\label{rem-rel-explor}
{\rm The above gives a highly-convenient {\em joint construction} of the first-passage percolation together with the graph $\CMnd$ itself. It is this exploration that was crucially used in \cite{BhaHofHoo12}. It is close to the one used in \cite{BhaHofHoo09b}, which also \cite{BarHofKom17} relies on, but there the cycle-removal procedure was not applied. Since we will use the exploration only up to $a_n=n^{\rho}$ vertices with $\rho>0$ small, in fact it is unlikely that any cycles appear up to this point (see e.g., \cite[Lemma 4.2]{Hofs24}).}\hfill $\ensymboldefinition$
\end{remark}

\begin{remark}[First-passage percolation exploration from two sources]
\label{rem-explor-two-sources}
{\rm On several occasions, it will be useful to explore from at least two, or even more, sources simultaneously, either with the same edge-weights, or independent edge-weights. While we do not write such extensions out in full detail, this can easily be done, for example, by performing the explorations from the different sources at the same time or one-by-one. When exploring from the two sources one-by-one, the above then describes the exploration from the first source, for that from another source, the only possible additional issue is that the edge-weights of edges {\em connecting} between the various SWGs need to be determined consistently. This is why in \cite{BhaHofHoo12}, the half-edges are again tested to see whether they go to new vertices, or to vertices that are already in one SWG, and the edge-weights are fixed accordingly. When exploring at the same time, we will write $\SWG^{\sss(1)}(t)$ and $\SWG^{\sss(2)}(t)$ for the sets of vertices that can be reached at time $t$ from the two sources $\Ver_1$ and $\Ver_2$, respectively, and $\SWG^{\sss(1,2)}(t)$ for their union.
}\hfill $\ensymboldefinition$
\end{remark}

\begin{remark}[Coupling to branching processes]
\label{rem-coupling}
{\rm The above construction also allows for a convenient coupling between the exploration in the random graph and that on the corresponding $n$-dependent branching process tree. Indeed, when pairing a half-edge, we draw a half-edge uar from the collection of {\em all} half-edges, and when it is still available, then we pair it in the same manner for tree and graph. This means that differences between graph and branching process exploration arise due to half-edges being redrawn for the branching process, as well as the cycle-removal procedure described above. This is worked out in full detail in \cite[Section 4.2]{Hofs24}.}\hfill $\ensymboldefinition$
\end{remark}

\subsection{The epidemic curve for first-passage percolation}
\label{sec-epidemic-curve}


We will see that at time $T_{a_n}+t$, and due to the explosive nature of our unimodular age-dependent branching process, a {\em positive proportion} of the vertices will be found by the first-passage percolation exploration. To describe how the first-passage percolation sweeps through the graph, we need some further notation. Let
	\eqn{
	\barNtk=\#\{v\colon d_v=k \mbox{ and } v\in \SWG(T_{a_n}+t)\}/n
	}
denote the fraction of vertices that have degree $k$ and that have been reached by the flow at time $T_{a_n}+t$. Further, for an edge $e=xy$ consisting of two half-edges $x$ and $y$ that are incident to vertices $v_x$ and $v_y$, we say that $e$ {\em spreads the flow at time $s$} when $v_x$ (or $v_y$) is occupied by first-passage percolation at time $s$, and $v_y$ (or $v_x$) is then occupied by first-passage percolation at time $s$ {\em through the edge $e$}. We let
	\eqn{
	\barLt=\#\{e\colon \mbox{ $e$ has spread the flow at time }T_{a_n}+t\}/[L_n/2]
	}
denote the proportion of edges that have spread the first-passage percolation exploration by time $T_{a_n}+t$. We further let $V(k)$ denote the explosion time of a vertex that has $k$ children in a first-passage percolation process on a discrete branching process with edge-weights distributed as $X$, so that the explosion time of the one-stage $(D^\star-1, X)$ branching process has distribution $V(D^\star-1)$, while that of a unimodular $(D,X)$ branching process has distribution $V(D)$. Our main result describing how first-passage percolation sweeps through the graph is as follows:

\begin{theorem}[First-passage percolation epidemic curve]
\label{thm-degree-time-FPP}
As $n\rightarrow \infty$,
	\eqn{
	\label{Ntk-conv-FPP}
	\barNtk\convp \prob(V(k)\leq t)\prob(D=k).
	}
Consequently, the proportion of vertices occupied by first-passage percolation at time $T_{a_n}+t$ converges in probability to $\prob(V(D)\leq t)$. Further,
	\eqn{
	\label{Lt-conv-FPP}
	\barLt \convp \prob\big(X+\tildeV_a\wedge\tildeV_b\leq t\big),
	}
where $(\tildeV_a,\tildeV_b)$ are two independent copies of $V(D^\star-1)$.
\end{theorem}

Theorem \ref{thm-degree-time-FPP} is interesting in its own right, as it implies the existence of an {\em epidemic curve}, as also derived in \cite{BhaHofKom14} in the setting of first-passage percolation on the configuration model with finite-variance degrees. Indeed, \eqref{Ntk-conv-FPP} together with the convergence $T_{a_n}\convd V$, where $V$ is the explosion time of the unimodular $(D,X)$ branching process, shows that the fraction of vertices occupied by first-passage percolation at time $t$ (or infected vertices when interpreting first-passage percolation as a model for disease spread) converges in distribution to $\prob(V^{\sss(2)}<t-V\mid V)=F_{\sss V}(t-V)$ for every $t\in \mathbb{R}$, where $F_{\sss V}$ is the distribution function of $V$. We refer to \cite{BhaHofKom14} for an extensive discussion of the consequences of such a result. A related version of \eqref{Ntk-conv-FPP} was proved in \cite[Proposition 29]{HofShn23}, but we crucially rely on the extension to fixed $k$ in \eqref{Ntk-conv-FPP}, as well as \eqref{Lt-conv-FPP}.

We extend Theorem \ref{thm-degree-time-FPP} to the winning type in Proposition \ref{prop-degree-time} below. Thus, the winning type sweeps through the graph in a similar way as first-passage percolation would {\em in the absence of competition}. The reason for this is that the winning type, instantly after having exploded, reaches a positive proportion of the vertices in the graph (and particularly the high-degree vertices), thus effectively blocking the spread of the losing type.


\subsection{First-passage percolation preliminaries}
\label{sec-prel}
In this section, we collect some results on first-passage percolation on scale-free configuration models, proved in earlier work. Our first preliminary result allows us to treat the degrees found in the random graph exploration as iid random variables, at least for a long time. In its statement, we let $(B_i^{\sss(n)})_{i\geq 1}$ denote the forward degrees of the vertices occupied by the first-passage exploration. That is, let $v_m$ denote the unique vertex that is in $\SWG(T_m)$ but not in $\SWG(T_{m-1})$ (and recall from Section \ref{sec-SWGs} that there are $m+1$ vertices in $\SWG(T_m)$, so the source of the first-passage percolation is not counted in $(v_m)_{m\geq 1}$), so that $B_m^{\sss(n)}=d_{v_m}-1$:

\begin{proposition}[Bhamidi, van der Hofstad, Hooghiemstra {\protect \cite{BhaHofHoo09b}}]
\label{prop-prel}
Consider first-passage percolation on a graph generated by the configuration model with iid degrees whose distribution satisfies \eqref{F-conditions}--\eqref{gc}. Then, there exists a $\rho>0$ such that the sequence $(B_i^{\sss(n)})_{i\geq 1}$ can be coupled to the iid sequence $(D^\star_i-1)_{i\geq 1}$ with law \eqref{forwarddegree} in such a way that $(B_i^{\sss(n)})_{i=1}^{n^\rho}=(D^\star_i-1)_{i=1}^{n^\rho}$ whp. Moreover, the vertices corresponding to these draws are whp distinct.\\
The above construction, with independent forward degrees, can be extended to the first-passage exploration from  two sources $\Ver_1$ and $\Ver_2$ chosen independently and uar from $[n]$.
\end{proposition}

\proof The coupling is proved in \cite[Proposition 4.7]{BhaHofHoo09b}. Note that our numbering is slightly different than in  \cite{BhaHofHoo09b}, since  for us there are $m+1$ vertices in $\SWG(T_m)$. The fact that the vertices are distinct is \cite[Lemma A.1]{BhaHofHoo09b}. The extension to two sources and the resulting independence follows from \cite[Proposition 4.8]{BhaHofHoo09b}.
\qed
\medskip

In our second result, we list some properties of the first-passage percolation exploration through the network. In its statement, we write $a_n=n^\rho$, where $\rho>0$ is taken from Proposition \ref{prop-prel}:

\begin{proposition}[Baroni, van der Hofstad, Komj\'athy \cite{BarHofKom17}]
\label{prop-prel-2}
Consider first-passage percolation on a graph generated by the configuration model with iid degrees whose distribution satisfies \eqref{F-conditions}--\eqref{gc}. Let $X^{\sss(1)}$ and $X^{\sss(2)}$ be such that at least one of the unimodular $(D,X^{\sss(1)})$ and $(D,X^{\sss(2)})$ branching processes is explosive. 
Let $m=m_n\to\infty$ with $m_n\leq a_n$, and fix two sources $\Ver_1$ and $\Ver_2$ chosen independently and uar from $[n]$. Then $(T_{m_n}^{\sss(1)}, T_{m_n}^{\sss(2)})\convd (V^{\sss(1)},V^{\sss(2)})$ as $n\to\infty$, where $(V^{\sss(1)},V^{\sss(2)})$ are independent explosion times of the unimodular $(D,X^{\sss(1)})$ and $(D,X^{\sss(2)})$ branching processes (with $V^{\sss(i)}=\infty$ when $(D,X^{\sss(i)})$ is conservative). Further, $\SWG^{\sss(1)}(T_{m_n}^{\sss(1)})$ and $\SWG^{\sss(2)}(T_{m_n}^{\sss(2)})$ are whp disjoint.
\end{proposition}

\proof This is a combination of \cite[Theorem 4]{BarHofKom17} and its proof, and the coupling in \cite[Proposition 4.7]{BhaHofHoo09b} discussed in Proposition \ref{prop-prel} above. For the independence, it is useful to note that $(X^{\sss(1)}_e)_{e\in E(\CMnd)}$ and $(X^{\sss(2)}_e)_{e\in E(\CMnd)}$ are independent, so that the asymptotic independence of the graph explorations in Proposition \ref{prop-prel} implies the independence of the explosion times.
\qed

\begin{remark}[A source with fixed degree]
\label{rem-source-degree-k}
{\rm Fix $k\geq 1$ such that $\prob(D=k)>0$. Then, the above result remains true when we start first-passage percolation from a vertex $\Ver_2$ chosen uar from all vertices of degree $k$. In that case, $V^{\sss(2)}$ is replaced by $V^{\sss(2)}(k)$, which is the explosion time of a unimodular $(D,X^{\sss(2)})$ process conditionally on $D=k$. This will be needed later on.
}\hfill$\ensymboldefinition$
\end{remark}



\subsection{The epidemic curve: Proof of Theorem \ref{thm-degree-time-FPP}}
\label{sec-proof-epidemic-curve}

In this section, we use Propositions \ref{prop-prel} and \ref{prop-prel-2} to prove Theorem \ref{thm-degree-time-FPP}. This proof is a minor modification of \cite[Proof of Proposition 4.2]{DeiHof16}, where a similar statement was proved for the winning type (see also Proposition \ref{prop-degree-time} below).

We start with the proof of \eqref{Ntk-conv-FPP}. Let $\Ver$ be a vertex chosen uar from $[n]$, independently of the first-passage percolation source $\Ver_1$, and write $\indicwo{\sss \Ver}^{\sss(t,k)}$ for the indicator taking the value 1 when vertex $\Ver$ has degree $k$ and was occupied by the first-passage percolation exploration from $\Ver_1$ at time $T_{a_n}+t$. Note that, with $G_n$ denoting the realisation of the configuration model including its edge-weights,
	\eqn{
	\label{Ntk-repr}
	\barNtk=\expec[\indicwo{\sss \Ver}^{\sss(t,k)}\mid G_n].
	}
We will show that $\E[\indicwo{\sss \Ver}^{\sss(t,k)}\mid G_n]\convp \prob(V(k)\leq t)\prob(D=k)$ by using a conditional second moment method. We perform the analysis conditionally on $\SWG(T_{a_n})$, that is, the first-passage percolation exploration graph at the time when it reaches size $a_n$ and $a_n$ is chosen as $n^{\rho}$ with $\rho>0$ as in Proposition \ref{prop-prel}. To apply the conditional second moment method, first note that
	\eqan{
	\expec[\barNtk \mid \SWG(T_{a_n})]&=\prob(\indicwo{\sss \Ver}^{\sss(t,k)}=1 \mid \SWG(T_{a_n}))\nn\\
	&=\prob(\mbox{$\Ver$ is found by flow at time $T_{a_n}+t$}
	\mid \SWG(T_{a_n}),
	D_{\sss \Ver}=k)\prob(D_{\sss \Ver}=k),
	}
and, with $\Ver'$ be chosen uar from $[n]$ and independently of $\Ver, \Ver_1$,
	\eqan{
	\expec[(\barNtk)^2 \mid \SWG(T_{a_n})]
	&=\prob(\indicwo{\sss \Ver}^{\sss(t,k)}=\indicwo{\sss \Ver'}^{\sss(t,k)}=1 \mid \SWG(T_{a_n}))\nn\\
	&=\prob(\mbox{$\Ver,\Ver'$ are found by flow at time $T_{a_n}+t$}\mid
	\SWG(T_{a_n}, D_{\sss \Ver}=D_{\sss \Ver'}=k)\nn\\
	&\qquad\times \prob(D_{\sss \Ver}=k)\prob(D_{\sss \Ver'}=k),
	}
where we use that the events $\{D_{\sss \Ver}=k\}$ and $\{D_{\sss \Ver'}=k\}$ are independent. Therefore, it suffices to show that the first factors in the above two rhs's converge to $\prob(V(k)\leq t)$ and $\prob(V(k)\leq t)^2$, respectively. Indeed, in this case,
	\eqn{
	\expec[\barNtk\mid \SWG(T_{a_n})]\convp \prob(V(k)\leq t) \prob(D=k),
	}
while ${\mathrm{Var}}(\barNtk\mid \SWG(T_{a_n}))=\op(1)$, so that $\barNtk\convp \prob(V(k)\leq t) \prob(D=k)$, as required.


Recall from Section \ref{sec-SWGs} that $W_n(\Ver, \Ver_1)$ denotes the passage time between vertices $\Ver_1$ and $\Ver$. It follows from the analysis in \cite{BhaHofHoo09b, BarHofKom17}, summarised in Propositions \ref{prop-prel}-\ref{prop-prel-2}, that $W_n(\Ver, \Ver_1)$ converges in distribution to $V+V(k)$ (recall Remark \ref{rem-source-degree-k}).  As described above, we first grow $\SWG(T_{a_n})$ from $\Ver_1$ as described in Section \ref{sec-SWGs}. Then we grow the SWG from $\Ver$ until it hits $\SWG(T_{a_n})$. This occurs at a time that converges in distribution to $V(k)$ -- indeed, $V(k)$ describes the asymptotic explosion time for an exploration process started at a vertex with degree $k$. Hence,
	\begin{equation}\label{eq:onetypeconv}
	\prob(W_n(\Ver, \Ver_1)\leq T_{a_n}+t\mid \SWG(T_{a_n}), D_{\sss U}=k)
	\convp \prob(V(k)\leq t).
	\end{equation}
In a similar way, we conclude that 
	\eqn{
	\prob(W_n(\Ver_1, \Ver),W_n(\Ver_1, \Ver')\leq T_{a_n}+t\mid \SWG(T_{a_n}), D_{\sss \Ver}=D_{\sss \Ver'}=k)\convp \prob(V(k)\leq t)^2.
	}
The second moment method now proves \eqref{Ntk-conv-FPP}.

The proof of \eqref{Lt-conv-FPP} is similar. Indeed, instead of \eqref{Ntk-repr}, we now start from $\barLt=\expec[\indicwo{e}^{\sss(t)}\mid G_n]$, where $e$ is a uniform edge in the graph and $\indicwo{e}^{\sss(t)}$ denotes the indicator that $e$ spreads the infection before time $t$. We then again use a conditional second moment, and note that 	
	\[
	\expec[\barLtwin \mid \SWG(T_{a_n})]
	=\prob(e\text{ has spread the flow by time }T_{a_n}+t \mid \SWG(T_{a_n})).
	\]
In this expectation, a uniform edge can be obtained by drawing a half-edge uar, and pairing it to a uniform other half-edge. Let $a$ and $b$ be the vertices at the two ends of $e$, and let $\widetilde V_a$ and $\widetilde V_b$ be the explosion times of the vertices $a$ and $b$, respectively, when the type 1 infection is not allowed to use the edge $e$. Then, $e$ has spread the infection by time $T_{a_n}+t$ precisely when either the explosion time $\widetilde V_a$ plus the edge-weight $X_e$ are at most $t$ (in which case, $a$ is first occupied by the flow and then spreads the flow to vertex $b$), or the explosion time $\widetilde V_b$ plus the edge-weight $X_e$ are at most $t$ (in which case, $b$ is first occupied by the flow and then spreads the flow to vertex $a$). We conclude that $\expec[\barLt \mid \SWG(T_{a_n})]\convp \prob\big(X+\widetilde V_a\wedge\widetilde V_b\leq t\big)$. The extensions to the second-moment computations are the same as for $\barNtk$.
\qed

\section{The winning type occupies most: Proof of Theorem \ref{thm-main}(a)}
\label{sec-most}
In this section, we return to the competition problem, where each edge has two edge-weights associated at them.
We complete the proof of Theorem \ref{thm-main}(a). This section is organised as follows. In Section \ref{sec-identifyinfg-winning-type}, we start by identifying the winning type and give the overview of proof. In Section \ref{sec-around-explosion-disjointness}, we show that $\SWG^{\sss(1)}(T_{a_n/2}^{\sss(1)}\wedge T_{a_n/2}^{\sss(2)}+\varepsilon)$ and $\SWG^{\sss(2)}(T_{a_n/2}^{\sss(1)}\wedge T_{a_n/2}^{\sss(2)}+\varepsilon)$ are whp disjoint, and that they are even at distance at least 2 (so that there are no direct edges between their vertices). 
Let $\mathcal{W}_n\in \{1,2\}$ denote the winning type. In Section \ref{sec-size-biased-degrees}, we investigate the degree distribution of the complement of $\SWG^{\sss(\mathcal{W}_n)}(T_{a_n/2}^{\sss(1)}\wedge T_{a_n/2}^{\sss(2)})$, and prove that it has all moments. Because of this, the losing type is seriously hindered in its exploration, and we show in Section \ref{sec-around-explosion-expensive} that even without the further competition with the winning type, the losing type restricted to the complement of $\SWG^{\sss(\mathcal{W}_n)}(T_{a_n/2}^{\sss(1)}\wedge T_{a_n/2}^{\sss(2)})$ does not explode after this time. In Section \ref{sec-proof-thm-a}, we then complete the proof of Theorem \ref{thm-main}(a). In Section \ref{sec-winning-hardly-hindered-losing}, we show that the winning type sweeps through the graph as if the losing type were not there, see Proposition \ref{prop-degree-time} and compare it to Theorem \ref{thm-degree-time-FPP}. 
This result will be crucial in Section \ref{sec-almost-all}, where we identify the limit in distribution of the number of vertices occupied by the losing type.

\subsection{Identifying the winning type and overview proof of Theorem \ref{thm-main}(a)}
\label{sec-identifyinfg-winning-type}
We start by identifying the winning type. Indeed, as explained in Section \ref{sec-overview}, $\SWG^{\sss(1,2)}(T_{a_n/2}^{\sss(1)}\wedge T_{a_n/2}^{\sss(2)})$ whp consists of two disjoint trees (recall Propositions \ref{prop-prel} and \ref{prop-prel-2}). One of these trees is whp large, the other is small. The larger one corresponds to the type $i$ for which $T_{a_n/2}^{\sss(i)}$ is the smallest. Thus, in this section, we condition on the event $A_n=\{T_{a_n/2}^{\sss(1)}<T_{a_n/2}^{\sss(2)}\}$. Recall that $\bar{N}_1(n)$ denotes the proportion of vertices occupied by type 1. To prove Theorem \ref{thm-main}(a), it then suffices to show that $\bar{N}_1(n)=1+\op(1)$ in this case, as we know that $(T_{a_n/2}^{\sss(1)},T_{a_n/2}^{\sss(2)})\convd (V^{\sss(1)}, V^{\sss(2)})$. 

The proof that $\bar{N}_1(n)=1+\op(1)$ uses a first moment method on $\bar{N}_1(n),$ conditionally on $A_n$. We start by investigating the conditional first moment, which can be written as
	\eqn{
	\expec\big[\bar{N}_1(n)\mid T_{a_n/2}^{\sss(1)}<T_{a_n/2}^{\sss(2)}\big]
	=\prob(\Ver \text{ \rm is type 1 infected} \mid T_{a_n/2}^{\sss(1)}<T_{a_n/2}^{\sss(2)}).
	} 
The main aim of this section will be to show that
	\eqn{
	\label{main-aim-sec2}
	\prob(\Ver \text{ \rm is type 1 infected} \mid T_{a_n/2}^{\sss(1)}<T_{a_n/2}^{\sss(2)})\rightarrow 1,
	}
where $\Ver$ is a vertex chosen uar from $[n]$, independently of $\Ver_1$.
Markov's inequality on $1-\bar{N}_1(n)$ then shows that in fact $\bar{N}_1(n)\convp 1$, as required. The remainder of this section is devoted to the proof of \eqref{main-aim-sec2}. It is divided into several propositions.

\subsection{Structure shortly after explosion: disjointness of trees}
\label{sec-around-explosion-disjointness}
We start by showing that $\SWG^{\sss(1)}(T_{a_n/2}^{\sss(1)}\wedge T_{a_n/2}^{\sss(2)}+\varepsilon)$ and $\SWG^{\sss(2)}(T_{a_n/2}^{\sss(1)}\wedge T_{a_n/2}^{\sss(2)}+\varepsilon)$ are whp disjoint when $\vep>0$ is small. We emphasise that the smallest-weight graphs discussed here all {\em ignore the competition}, we later include the competition between types 1 and 2 (see Remark \ref{rem-comp-vs-fpp}):

\begin{proposition}[SWGs remain disjoint for a small time]
\label{prop-disj-trees}
As $\vep\searrow 0$,
	\be\label{emptin}
	\limsup_{n\rightarrow\infty}\prob(\SWG^{\sss(1)}(T_{a_n/2}^{\sss(1)}\wedge T_{a_n/2}^{\sss(2)}+\varepsilon)\cap \SWG^{\sss(2)}(T_{a_n/2}^{\sss(1)}\wedge T_{a_n/2}^{\sss(2)}+\varepsilon)=\varnothing)=1-o(1),
	\ee
where $a_{n}=n^{\rho}$ with $\rho$ identified in Proposition \ref{prop-prel}. Further, it is unlikely that there exists a direct edge between $\SWG^{\sss(1)}(T_{a_n/2}^{\sss(1)}\wedge T_{a_n/2}^{\sss(2)}+\varepsilon)$ and $\SWG^{\sss(2)}(T_{a_n/2}^{\sss(1)}\wedge T_{a_n/2}^{\sss(2)}+\varepsilon)$ when $\vep>0$ is small.
\end{proposition}

\begin{proof} We rely on Theorem \ref{thm-degree-time-FPP}, which implies that, on the event $\{T_{a_n/2}^{\sss(1)}<T_{a_n/2}^{\sss(2)}\}$,
	\eqn{
	\frac{1}{n}|\SWG^{\sss(1)}(T_{a_n/2}^{\sss(1)}\wedge T_{a_n/2}^{\sss(2)}+\varepsilon)|\convp  \prob(V^{\sss(1)}(D)\leq \vep).
	}
Further, for $\vep>0$ sufficiently small, whp for $\vep>0$ small, $\SWG^{\sss(2)}(T_{a_n/2}^{\sss(1)}\wedge T_{a_n/2}^{\sss(2)}+\varepsilon)=\SWG^{\sss(2)}(T_{a_n/2}^{\sss(1)}\wedge T_{a_n/2}^{\sss(2)})$, where $|\SWG^{\sss(2)}(T_{a_n/2}^{\sss(1)}\wedge T_{a_n/2}^{\sss(2)})|$ is a tight random variable on the event that $T_{a_n/2}^{\sss(1)}<T_{a_n/2}^{\sss(2)}$. Since $ \prob(V^{\sss(1)}(D)\leq \vep)=o(1)$ as $\vep\searrow 0$, the fraction of vertices in $\SWG^{\sss(1)}(T_{a_n/2}^{\sss(1)}\wedge T_{a_n/2}^{\sss(2)}+\varepsilon)$ is small. The vertices in $\SWG^{\sss(2)}(T_{a_n/2}^{\sss(1)}\wedge T_{a_n/2}^{\sss(2)}+\varepsilon)$ consist of $\Ver_2$ together with a {\em tight} number of vertices whose degrees converge in distribution (recall Proposition \ref{prop-prel}). Given its degree, any vertex of that degree is equally likely to be in $\SWG^{\sss(2)}(T_{a_n/2}^{\sss(1)}\wedge T_{a_n/2}^{\sss(2)}+\varepsilon)$, and this is independent of whether it is in $\SWG^{\sss(1)}(T_{a_n/2}^{\sss(1)}\wedge T_{a_n/2}^{\sss(2)}+\varepsilon)$, since $(X^{\sss(1)}_e)_{e}$ and $(X^{\sss(2)}_e)_{e}$ are independent. The probability that any of the tight number of vertices in $\SWG^{\sss(2)}(T_{a_n/2}^{\sss(1)}\wedge T_{a_n/2}^{\sss(2)})$ is also in $\SWG^{\sss(1)}(T_{a_n/2}^{\sss(1)}\wedge T_{a_n/2}^{\sss(2)}+\varepsilon)$ is thus small for $\vep>0$ small. 

The proof that whp there are no edges linking $\SWG^{\sss(1)}(T_{a_n/2}^{\sss(1)}\wedge T_{a_n/2}^{\sss(2)}+\varepsilon)$ and $\SWG^{\sss(2)}(T_{a_n/2}^{\sss(1)}\wedge T_{a_n/2}^{\sss(2)}+\varepsilon)$ is identical, since the number of half-edges incident to $\SWG^{\sss(2)}(T_{a_n/2}^{\sss(1)}\wedge T_{a_n/2}^{\sss(2)}+\varepsilon)$ is also tight on the event $\{T_{a_n/2}^{\sss(1)}<T_{a_n/2}^{\sss(2)}\}$.
\end{proof}

\begin{remark}[Competition versus first-passage percolation]
\label{rem-comp-vs-fpp}
{\rm Since it is unlikely that there exists a direct edge between $\SWG^{\sss(1)}(T_{a_n/2}^{\sss(1)}\wedge T_{a_n/2}^{\sss(2)}+\varepsilon)$ and $\SWG^{\sss(2)}(T_{a_n/2}^{\sss(1)}\wedge T_{a_n/2}^{\sss(2)}+\varepsilon)$ when $\vep>0$ is small, the presence of competition is not noticed by the exploration of the graph by the two types. As such, the exploration from vertices $\Ver_1$ and $\Ver_2$ according to the edge-weights $(X^{\sss(1)}_e)_{e}$ and $(X^{\sss(2)}_e)_{e}$ is the {\em same} up to time $T_{a_n/2}^{\sss(1)}\wedge T_{a_n/2}^{\sss(2)}+\varepsilon$ as that of first-passage percolation from these sources with these edge-weights. Thus, up to time $T_{a_n/2}^{\sss(1)}\wedge T_{a_n/2}^{\sss(2)}+\varepsilon$, we can effectively use the first-passage percolation results in Propositions \ref{prop-prel} and \ref{prop-prel-2}.
}\hfill$\ensymboldefinition$
\end{remark}

\subsection{Structure of outside of winning SWG shortly after explosion: Size-biased reordering according to degrees}
\label{sec-size-biased-degrees}
Before moving on to our second main result, we investigate the degree distribution of the vertices that have not been found by first-passage percolation from $U_1$ with edge-weights $(X^{\sss(1)}_e)_{e}$. Let  $\tilde{\boldsymbol{d}}(\vep)$ be the (random) degree distribution of the vertices in $[n]\setminus \SWG^{\sss(1)}(T_{a_n/2}^{\sss(1)}+\varepsilon)$, and let $\tilde D_{n}(\vep)$ denote the degree of a uniformly chosen vertex in that set. Let $\field_{a_n}$ denote the $\sigma$-field generated by $\SWG^{\sss(1)}(T_{a_n/2}^{\sss(1)}+\varepsilon)$ and the full exploration process up to time $T_{a_n/2}^{\sss(1)}+\varepsilon$. The following proposition shows that $\tilde{\boldsymbol{d}}(\vep)$ has a tight second moment:

\begin{proposition}[Complement of SWG has bounded second moment degrees]
\label{prop-bd-sec-mom}
For all $\varepsilon>0$, $\tilde{\boldsymbol{d}}(\vep)$ satisfies that $\expec[\tilde D_{n}(\vep)^{2}\mid \field_{a_n}]$ is a tight sequence of random variables. 
\end{proposition}

\begin{proof} Recall Remark \ref{rem-comp-vs-fpp}, which implies that up to time $T_{a_n/2}^{\sss(1)}\wedge T_{a_n/2}^{\sss(2)}+\varepsilon$, there is whp no difference between the first-passage percolation and competition processes. We can thus rely on first-passage percolation estimates. We need to prove that, for every $\eta>0$, there exists a $K>0$ such that
	\eqn{
	\label{aim-tightness-second-moment}
	\prob\Big(\expec[\tilde D_{n}(\vep)^{2}\mid \field_{a_n}]>K\Big)\leq \eta.
	}
We note that, at time $T_{a_n/2}^{\sss(1)}+\varepsilon$, by \eqref{Ntk-conv-FPP} in Theorem \ref{thm-degree-time-FPP}, $\SWG^{\sss(1)}(T_{a_n/2}^{\sss(1)}+\varepsilon)$ has size
	\eqn{
	\label{Nn-vep-asympt}
	N_n^{\sss(\vep)}=n\sum_{k\geq 2} \bar{N}_n^{\sss (k,\vep)}=n\prob(V^{\sss(1)}(D)\leq \vep)(1+\op(1)).
	}
We conclude that whp, the event ${\mathcal E}_n =\{N_n^{\sss(\vep)}\in n[\prob(V^{\sss(1)}(D)\leq \vep)/2,2\prob(V^{\sss(1)}(D)\leq \vep)]\}\cap \{L_n\leq 2n \expec[D]\}$ holds. We then prove that whp
	\eqn{
	\label{aim-tightness-second-moment-rep}
	\prob_n\Big(\indicwo{{\mathcal E}_n}\expec_n[\tilde D_{n}(\vep)^{2}\mid \field_{a_n}]>K\Big)\leq \eta,
	}
where we let $\prob_n$ denote the law of the configuration model {\em given} the iid degrees $\boldsymbol{d}=(d_1,\dots, d_n)$, and $\expec_n$ the expectation wrt $\prob_n$.

In order to analyse $\tilde D_{n}(\vep)$, we introduce the following notation:

\begin{definition}[Size-biased reordering] 
\label{def-size-biased-reor}
{\rm Let $(d_v)_{v\in[n]}$ be a sequence of integers. The {\em uniform size-biased reordering} of $[n]$ wrt $(d_v)_{v\in[n]}$ is the random permutation $(\pi_i)_{i\in[n]}$ of $[n]$, defined recursively by
	\eqn{
	\label{size-biased-reor-1}
	\prob(\pi_1=v)=\frac{1}{n},
	}
and, given $\pi_{[i-1]}\equiv (\pi_1, \ldots, \pi_{i-1})$ and all $v\not\in {\sf SB}(i-1)\equiv \{\pi_1, \ldots, \pi_{i-1}\}$,
	\eqn{
	\label{size-biased-reor-i}
	\prob(\pi_i=v\mid \pi_{[i-1]})=\frac{d_v}{\sum_{u\in [n]\setminus {\sf SB}(i-1)}d_u}.
	}
}\hfill $\ensymboldefinition$
\end{definition}

The crux of Definition \ref{def-size-biased-reor} is that vertices keep being added with a probability proportional to their degrees, but without replacement. This makes the probability that vertices are added larger and larger when they have escaped being drawn for a long time. Further, due to the {\em size-biasing} effect in \eqref{size-biased-reor-i} wrt the degrees $(d_v)_{v\in[n]}$, vertices with higher degree are more likely to be chosen {\em earlier} in the exploration. 

In the usual size-biased reordering, \eqref{size-biased-reor-1} is replaced with $\prob(\pi_1=v)=d_v/\sum_{u\in[n]}d_u$, which would correspond to also $\pi_1$ being chosen in a size-biased manner. Definition \ref{def-size-biased-reor} has been set up in such a way that the vertices chosen in the exploration of a random vertex $\Ver_1$, and removing repetitions due to cycles, follow exactly these rules (recall the first-passage percolation exploration described in Section \ref{sec-SWGs}).

We conclude that the vertices in $\SWG^{\sss(1)}(T_{a_n/2}^{\sss(1)}+\varepsilon)$ form a collection of $N_n^{\sss(\vep)}$ uniform size-biased reordered vertices as in Definition \ref{def-size-biased-reor}. Indeed, each time a new vertex is being added, since the edge-weights are attached to the half-edge that already is in $\SWG^{\sss(1)}(T_{a_n/2}^{\sss(1)}+\varepsilon)$, the half-edge is chosen uar, conditionally on being incident to a vertex not yet in the SWG, and such a collection is precisely a uniform size-biased reordered vector. Further, its size is precisely $N_n^{\sss(\vep)}$ by construction.

Recall that $\tilde D_{n}(\vep)$ denotes the degree of a vertex chosen uar from $[n]\setminus \SWG^{\sss(1)}(T_{a_n/2}^{\sss(1)}+\varepsilon)$. The law of $\tilde D_{n}(\vep)$ conditionally on $\field_{a_n}$ is the same as the law of $\tilde D_{n}(\vep)$ conditionally on ${\sf SB}(N_n^{\sss(\vep)})$. As a result, 
	\[
	\expec_n[\tilde D_{n}(\vep)^{2}\mid \field_{a_n}]
	=\expec_n[\tilde D_{n}(\vep)^{2}\mid {\sf SB}(N_n^{\sss(\vep)})]
	\] 
is a tight sequence of random variables when we can show that the probability that a vertex of relatively high degree is not in ${\sf SB}(N_n^{\sss(\vep)})$ is small. We then bound using Markov's inequality,
	\eqan{
	\label{tightness-second-moment-a}
	\prob_n\Big(\indicwo{{\mathcal E}_n}\expec_n[\tilde D_{n}(\vep)^{2}\mid \field_{a_n}]>K\Big)
	&\leq 
	\frac{1}{K}\expec_n\Big[\indicwo{{\mathcal E}_n}\expec_n[\tilde D_{n}(\vep)^{2}\mid \field_{a_n}]\Big]\\
	&=\frac{1}{K}\expec_n\big[\indicwo{{\mathcal E}_n}\tilde D_{n}(\vep)^2\big]\nn\\
	&=\frac{1}{K}\expec_n\Big[\indicwo{{\mathcal E}_n}\expec_n[\tilde D_{n}(\vep)^2\mid N_n^{\sss(\vep)}\big]\Big],\nn
	}
since ${\mathcal E}_n$ is measurable wrt $\field_{a_n}$, as well as $N_n^{\sss(\vep)}$ and $(d_v)_{v\in[n]}$.

Recall \eqref{Nn-vep-asympt}. 
Fix $\delta\in (0,1-\vep)$. We can bound
	\eqn{
	\expec_n[\tilde D_{n}(\vep)^2\mid N_n^{\sss(\vep)}=\lceil \delta n\rceil]
	= 
	\sum_{v} d_v^2 \expec_n\Big[\frac{\indic{v\not\in {\sf SB}(\delta n)}}{\sum_{u\not\in {\sf SB}(\delta n)} d_u}\Big]
	\leq \frac{1}{n-n\delta} \sum_{v} d_v^2 \prob_n(v\not\in {\sf SB}(\delta n)),
	}
since $d_v\geq 1$, so that
	\eqn{
	\sum_{u\not\in {\sf SB}(\delta n)} d_u\geq |{\sf SB}(\delta n)|\geq n-n\delta.
	}
Note that
	\eqn{
	\prob_n(v\not\in {\sf SB}(k)\mid v\not\in {\sf SB}(k-1))\leq 1-\frac{d_v}{L_n}.
	}
Thus, 
	\eqn{
	\prob_n(v\not\in {\sf SB}(\delta n))\leq \Big(1-\frac{d_v}{L_n}\Big)^{\delta n}\leq \e^{-d_v \delta n/L_n}.
	}
We thus arrive at
	\eqn{
	\expec_n[[\tilde D_{n}(\vep)^2\mid N_n^{\sss(\vep)}=\lceil \delta n\rceil]\leq \frac{1}{1-\delta} 
	\expec_n\Big[D_n^2 \e^{-D_n \delta n/L_n}\mid N_n^{\sss(\vep)}=\lceil \delta n\rceil\Big].
	}
As a result, since $L_n\leq 2n\expec[D]$ and $N_n^{\sss(\vep)}\leq 2\prob(V^{\sss(1)}(D)\leq \vep)n$ on ${\mathcal E}_n$, and with $\delta=2\prob(V^{\sss(1)}(D)\leq \vep)$,
	\eqan{
	\prob_n\Big(\indicwo{{\mathcal E}_n}\expec_n[\tilde D_{n}(\vep)^{2}\mid \field_{a_n}]>K\Big)
	&\leq \frac{1}{K(1-\delta)}\expec_n\Big[\indicwo{{\mathcal E}_n}\expec_n\Big[D_n^2 \e^{-D_n \delta n/L_n}\mid N_n^{\sss(\vep)}\Big]\Big]\\
	&\leq \frac{1}{K(1-\delta)}\expec_n\Big[\indicwo{{\mathcal E}_n}\expec_n\Big[D_n^2 \e^{-D_n \delta /(2\expec[D])}\mid N_n^{\sss(\vep)}\Big]\Big]\nn\\
	&\leq \frac{1}{K(1-\delta)}\expec_n\big[D_n^2 \e^{-D_n \delta /(2\expec[D])}\big].\nn
	}
This has finite expectation when taking the expectation wrt the iid degrees (and in fact, has any moment), which implies the claim in \eqref{aim-tightness-second-moment-rep}.
\end{proof}


\subsection{Structure shortly after explosion: expensive to avoid large SWG}
\label{sec-around-explosion-expensive}

We next prove that avoiding the SWG of vertex $\Ver_1$ is expensive. Let $U$ be a uniform vertex in $[n]\setminus \SWG^{\sss(1)}(T_{a_n/2}^{\sss(1)}+\varepsilon)$. We write $W^{\sss(2)}_n(\Ver_2,U)$ for the passage time from $\Ver_2$ to $U$ using the edge-weights $(X^{\sss(2)}_e)_e$:

\begin{proposition}[Expensive to avoid large SWG]
\label{prop-exp-avoid-SWG}
There exists $b_{n}\rightarrow \infty$ such that
	\be
	\prob\Big(W_{n}^{\sss (2)}(\Ver_2,U) \text{ in } [n]\setminus \SWG^{\sss(1)}(T_{a_n/2}^{\sss(1)}+\varepsilon)\leq b_{n} \mid T_{a_n/2}^{\sss(1)}<T_{a_n/2}^{\sss(2)}\Big)\rightarrow 0.
	\ee
\end{proposition}

We complete the proof of Proposition \ref{prop-exp-avoid-SWG} at the end of this section, and first explain the strategy.
We consider the graph $\CMnd\setminus\SWG^{\sss(1)}(T_{a_n/2}^{\sss(1)}+\varepsilon)$ and a uniformly chosen vertex $U$ in it. To prove Proposition \ref{prop-exp-avoid-SWG}, we aim to show that, in $\CMnd\setminus \SWG^{\sss(1)}(T_{a_n/2}^{\sss(1)}+\varepsilon)$,
	\be
	\label{weight-restricted-infty}
	W_{n}^{\sss(2)}(\Ver_2,\Ver)\convp \infty.
	\ee

Since, by Proposition \ref{prop-bd-sec-mom}, the degree distribution outside of $\SWG^{\sss(1)}(T_{a_n/2}^{\sss(1)}+\varepsilon)$ has finite second moment, we can couple the growth of $\SWG^{\sss(2)}(T_{a_n/2}^{\sss(1)}\wedge T_{a_n/2}^{\sss(2)}+t)$ in $[n]\setminus \SWG^{\sss(1)}(T_{a_n/2}^{\sss(1)}\wedge T_{a_n/2}^{\sss(2)}+\varepsilon)$ to an age-dependent $n$-dependent branching process. The following lemma implies that $\SWG^{\sss(2)}(T_{a_n/2}^{\sss(1)}\wedge T_{a_n/2}^{\sss(2)}+t)$ grows at most exponentially:

\begin{lemma}[Bounding the growth of the losing process] 
\label{lem-bd-growth-losing}
There exist tight random variables $\lambda_n(\vep), K_n(\vep)\in (0,\infty)$ such that, for all $t>0$,
	\eqn{
	\indic{T_{a_n/2}^{\sss(1)}<T_{a_n/2}^{\sss(2)}}
	\expec_n\Big[|\SWG^{\sss(2)}(T_{a_n/2}^{\sss(1)}\wedge T_{a_n/2}^{\sss(2)}+t)| \mid \field_{a_n}\Big]\leq K_n(\vep)\e^{\lambda_n(\vep) t}.
	}
\end{lemma}

\begin{proof} Let $T_{a_n/2}^{\sss(1)}<T_{a_n/2}^{\sss(2)}$. Let $S_0'(n)$ denote the number of active half-edges for the exploration from $\Ver_2$ (recall Section \ref{sec-SWGs}), which is a tight random variable. The random variable $|\SWG^{\sss(2)}(T_{a_n/2}^{\sss(1)}\wedge T_{a_n/2}^{\sss(2)}+t)|$ is bounded from above by $|\SWG(t)|$, which is a first-passage exploration started from $D_n'=\sum_{j=1}^{S_0'(n)} B^{\sss(\vep)}_j$ half-edges, and $(B^{\sss(\vep)}_j+1)_{j\geq 1}$ form a collection of size-biased reordered random variables from $[n]\setminus \SWG^{\sss(1,2)}(T_{a_n/2}^{\sss(1)}\wedge T_{a_n/2}^{\sss(2)}+\varepsilon)$ (recall Definition \ref{def-size-biased-reor}). 

This can be realised by noting that the $S_0'(n)$ half-edges each have a certain residual life-time, and the exploration is bounded by the setting where all residual life-times are replaced by 0. This boils down to immediately pairing these half-edges, and starting novel exploration processes from the sibling half-edges of the half-edges that are being paired to. We let the total number of unpaired half-edges after this construction be equal to $\tilde\ell_n(\vep)$, which is whp bounded from below by $L_n(1-\delta)$ for some $\delta>0$ that can be made small for small $\vep>0$. By construction, the vertices with their unpaired half-edges form a smaller configuration model with fewer vertices and edges.

For this smaller configuration model, we use {\em path-counting techniques}, as pioneered by Janson in \cite[Proof of Lemma 5.1]{Jans09b}, and follow the ideas in \cite[Section 5.3]{BhaHofHoo09b}. We follow \cite[Proof of Proposition 2.2(b)]{BhaHofHoo09b}, adapting the argument where needed. Let $(\BPn'(s))_{s\geq 0}$ denote the age-dependent branching process starting from $D_n'$ individuals, and having offspring distribution equal to $\tilde D_{n}^\star(\vep)-1$, the (size-biased distribution of $\tilde D_{n}(\vep)$) minus 1. Let
	\eqn{
	\YBP(t)=|\{v\colon v\in \BPn'(s) \text{ for some }s\leq t\}|,
	}
denote the total number of individuals ever born into the $\BPn'$ before time $t$, and let
	\eqn{
	\label{YSWT-def}
	\YSWG(t)=|\{v\colon v\in \Alive(s) \text{ for some }s\leq t\}|,
	}
denote the number of half-edges in the SWG that have ever been alive before time $t$. Here we recall from Section \ref{sec-SWGs} that ${\sf A}(t)$ denotes the set of active half-edges, which here is restricted to the vertices in $[n]\setminus \SWG^{\sss(1,2)}(T_{a_n/2}^{\sss(1)}\wedge T_{a_n/2}^{\sss(2)}+\varepsilon)$. Further, we emphasise the fact that all our explorations are performed in $[n]\setminus \SWG^{\sss(1,2)}(T_{a_n/2}^{\sss(1)}\wedge T_{a_n/2}^{\sss(2)}+\varepsilon)$, so outside $\SWG^{\sss(1,2)}(T_{a_n/2}^{\sss(1)}\wedge T_{a_n/2}^{\sss(2)}+\varepsilon)$.

The expected number of descendants in generation $k$ of a branching process starting with $\ell$ individuals equals $\ell\bar{\nu}_n(\vep)^{k}$, where 
	\eqn{
	\bar{\nu}_n(\vep):=\expec_n'[\tilde D_{n}^\star(\vep)-1]
	}
denotes the mean offspring, and $\expec_n'$ denotes the expectation conditionally on $\field_{a_n}$, the degrees $(d_v)_{v\in[n]}$ and $\SWG^{\sss(2)}(T_{a_n/2}^{\sss(1)}\wedge T_{a_n/2}^{\sss(2)})$. Here, we deal with a unimodular branching process $\BPn'$ where in the first generation the mean number of offspring equals $\expec_n'[D_n']$. $\GF^{\sss \star k}(t)-\GF^{\sss \star (k+1)}(t)$ equals the probability that an individual of generation $k$ is alive at time $t$. Together this yields
	\eqan{
	\expec_n'[|\BPn'(t)|]&=\sum_{k=1}^{\infty} \expec_n'[D_n']\bar{\nu}_n(\vep)^{k-1} [\GF^{\sss \star k}(t)-\GF^{\sss \star (k+1)}(t)],\\
	\expec_n'[\YBP(t)]&=2\expec_n'[D_n']+2\expec_n'[D_n']\sum_{k=1}^{\infty} \bar{\nu}_n(\vep)^{k-1} \GF^{\sss \star k}(t).
	}
To bound $\expec_n'[|\BPn'(t)|]$ and $\expec_n'[\YBP(t)]$ for various choices of $t$, according to \cite[(4.41)]{BhaHofHoo12b},
	\eqn{
	\label{expec-BPt}
	\expec[|\BPn'(t)|]\leq \expec_n'[D_n'] K_n(\vep)\e^{\lambda_n(\vep) t}(1+o(1)),
	}
where $\lambda_n(\vep)$ and $K_n(\vep)$ are functions of $\bar{\nu}_n(\vep)$. Our case involves settings where $\bar{\nu}_n(\vep)$ is uniformly bounded (it is a tight random variable by Proposition \ref{prop-bd-sec-mom}). This then yields an excellent upper bound for $\expec_n'[|\BPn'(t_n)|]$. A bound of the same order of magnitude for $\expec_n'[\YBP(t)]$ then straightforwardly follows from \cite[(5.30) in Lemma 5.4]{BhaHofHoo12b} that relates $\expec_n'[\YBP(t)]$ to $\expec_n'[|\BPn'(t)|]$

We continue with $\expec_n'[\YSWG(t)]$. We use the same steps as above,
and start by computing
	\eqn{
	\expec_n'[\YSWG(t)]=2\expec_n'[D_n']+2\expec_n'[D_n']\sum_{k=1}^{\infty} \GF^{\sss \star k}(t) \expec_n'[P_k^{\star}],
	}
where
	$$
	P_k^{\star}=\sum_{|\vec\pi|=k} (d_{\pi_k}-1),\qquad k\ge1,
	$$ 
is the sum of the number of  half-edges at the ends of paths $\vec\pi=(\pi_0, \ldots, \pi_k)$ of length $k$ in $\CMnd$, from one of the $D_n'$ initial half-edges, and restricted to the vertices in $[n]\setminus \SWG^{\sss(1,2)}(T_{a_n/2}^{\sss(1)}\wedge T_{a_n/2}^{\sss(2)}+\varepsilon)$ for $\pi_i$ with $i\in[k]$, while $\pi_0\in  \SWG^{\sss(2)}(T_{a_n/2}^{\sss(1)}\wedge T_{a_n/2}^{\sss(2)}+\varepsilon)$. Following \cite[Proof of Lemma 5.1]{Jans09b}, we find that
	\eqn{
	\label{res-janson}
	\expec_n'[P_k^{\star}]=\sum_{v_1, \ldots, v_k} \expec_n'[D'_n]\prod_{i=1}^k \frac{d_{v_i}(d_{v_i}-1)}{\tilde\ell_n(\vep)-2i+1}
	\leq \expec_n'[D'_n]\bar{\nu}_n(\vep)^{k},
	}
where the sum is taken over distinct vertices in $[n]\setminus \SWG^{\sss(1,2)}(T_{a_n/2}^{\sss(1)}+\varepsilon)$, and we recall that $\tilde\ell_n(\vep)$ denotes the number of unpaired half-edges in $[n]\setminus \SWG^{\sss(1,2)}(T_{a_n/2}^{\sss(1)}+\varepsilon)$.  We obtain, for $\bar{\nu}_n(\vep)\ge1$ (which we may assume without loss of generality),
	\eqn{
	\label{expec-YSWT}
	\expec_n'[\YSWG(t)]\leq 2\expec_n'[D_n']+2 \expec_n'[D'_n]\sum_{k=0}^{\infty} \GF^{\sss \star k}(t)\bar{\nu}_n(\vep)^k
	\le \expec_n'[\YBP(t)],
	}
so the requested bound follows from our bounds for $\expec_n'[\YBP(t)]$.
\end{proof}

\noindent
Now we are ready to complete the proof of Proposition \ref{prop-exp-avoid-SWG}:\\

\noindent
{\it Proof of Proposition \ref{prop-exp-avoid-SWG}.} We let $\prob_n'$ denote the conditional probability measure given $\field_{a_n}$, the degrees $(d_v)_{v\in[n]}$ and $\SWG^{\sss(2)}(T_{a_n/2}^{\sss(1)}\wedge T_{a_n/2}^{\sss(2)})$. We work on the event that $T_{a_n/2}^{\sss(1)}<T_{a_n/2}^{\sss(2)}$. Let $M_n^{\sss(\vep)}=|\SWG^{\sss(1,2)}(T_{a_n/2}^{\sss(1)}+\varepsilon)|$.

We compute that 
	\eqan{
	&\prob_n'\Big(W_{n}^{\sss (2)}(\Ver_2,U) \text{ in } [n]\setminus \SWG^{\sss(1,2)}(T_{a_n/2}^{\sss(1)}+\varepsilon)\leq b_{n}\Big)\nn\\
	&\qquad =\frac{1}{n-M_n^{\sss(\vep)}}\sum_{u\in [n]\setminus \SWG^{\sss(1,2)}(T_{a_n/2}^{\sss(1)}+\varepsilon)}
	\prob_n'\Big(W_{n}^{\sss (2)}(\Ver_2,u) \text{ in } [n]\setminus \SWG^{\sss(1,2)}(T_{a_n/2}^{\sss(1)}+\varepsilon)\leq b_{n}\Big).\nn\\
	&\qquad \leq \frac{1}{n-M_n^{\sss(\vep)}}\expec_n'\Big[|\SWG^{\sss(2)}(T_{a_n/2}^{\sss(1)}+b_n)|\Big].
	}
By Lemma \ref{lem-bd-growth-losing}, the expectation is bounded by $K_n(\vep)\e^{\lambda_n(\vep) b_n}$, while, by Theorem \ref{thm-degree-time-FPP}, $n-M_n^{\sss(\vep)}=\Thetap(n)$. We conclude that the rhs is $\op(1)$ when $b_n\leq \delta \log{n}$ for some $\delta$ sufficiently small, on the event that $\lambda_n(\vep)<1/\delta$, which occurs whp since $\lambda_n(\vep)$ and $K_n(\vep)$ are tight random variables.
\qed
\medskip

We close this section with a statement about the tightness of $|\SWG^{\sss(2)}(T_{a_n/2}^{\sss(1)}+b)|$ for all $b>0$:

\begin{corollary}[Losing type gets at most tight number of vertices]
\label{cor-losing-tight}
On the event that $T_{a_n/2}^{\sss(1)}<T_{a_n/2}^{\sss(2)}$, the random variable $|\SWG^{\sss(2)}(T_{a_n/2}^{\sss(1)}+b)|$ is tight for all $b>0$.
\end{corollary}

\proof This follows directly from Lemma \ref{lem-bd-growth-losing}, which proves that the expectation is whp bounded by $K_n(\vep)\e^{\lambda_n(\vep) b},$
which implies that $|\SWG^{\sss(2)}(T_{a_n/2}^{\sss(1)}+b)|$ is tight.
\qed

\subsection{Proof of Theorem \ref{thm-main}(a)} 
\label{sec-proof-thm-a}

We again condition on $T_{a_n/2}^{\sss(1)}<T_{a_n/2}^{\sss(2)}$. We apply a first moment method on $N_n^{\sss(2)}$, the number of vertices occupied by type 2. Recall \eqref{main-aim-sec2}, which shows that it suffices to prove that
	\eqn{
	\prob(\Ver \text{ \rm is type 2 infected} \mid T_{a_n/2}^{\sss(1)}<T_{a_n/2}^{\sss(2)})\rightarrow 0.
	}
We note that
	\eqn{
	\prob(\Ver \in \SWG^{\sss(2)}(T_{a_n/2}^{\sss(1)}+\varepsilon) \mid T_{a_n/2}^{\sss(1)}<T_{a_n/2}^{\sss(2)})\rightarrow 0,
	}
since $|\SWG^{\sss(2)}(T_{a_n/2}^{\sss(1)}+\varepsilon)|$ is whp tight for $\vep>0$ sufficiently small (recall Corollary \ref{cor-losing-tight}). Thus, $\Ver$ is unlikely to be occupied by type 2 before time $T_{a_n/2}^{\sss(1)}+\varepsilon$. Further, by Proposition \ref{prop-exp-avoid-SWG}, there exists a $b_n\ra \infty$ such that
	\eqn{
	\prob(\Ver \in \SWG^{\sss(2)}(T_{a_n/2}^{\sss(1)}+b_n) \mid T_{a_n/2}^{\sss(1)}<T_{a_n/2}^{\sss(2)})\rightarrow 0,
	}
Thus, $\Ver$ is unlikely to be occupied by type 2 before time $T_{a_n/2}^{\sss(1)}+b_n.$ By Theorem \ref{thm-degree-time-FPP}, whp $\Ver$ is discovered by the first-passage percolation flow from vertex $\Ver_1$ before time $b$ when $b$ is sufficiently large. However, it could be that this discovery is {\em unlucky} in that it is {\em blocked} by the competition from vertex $\Ver_2$. We next prove that this is unlikely to happen.

Consider the SWG from vertex $U$, where \whp $U\not\in\SWG^{\sss(2)}(\Toan+\vep)$ since $|\SWG^{\sss(2)}(\Toan+\vep)|$ is a tight random variable (again recall Corollary \ref{cor-losing-tight}). When $U\in \SWG^{\sss(1)}(\Toan+\vep)$, which, by Theorem \ref{thm-degree-time-FPP}, occurs with conditional probability $\prob(V^{\sss(1)}(k)\leq \vep)+\op(1)$, then it cannot be blocked by the fact that $\SWG^{\sss(1)}(\Toan+\vep)$ and $\SWG^{\sss(2)}(\Toan+\vep)$ are whp the same as their first-passage percolation counterparts that is obtained by ignoring the competition (recall Proposition \ref{prop-prel-2}). We thus focus on the case where $U\not\in\SWG^{\sss(1)}(\Toan+\vep)$. 

Without the presence of the type 2 infection, the exploration from such a $U$ hits $\SWG^{\sss(1)}(\Toan+\vep)$ when it has reached some random size $Q_n$, where $Q_n$ is a tight random variable for every $\vep>0$. 
Indeed, recall the exploration process defined in Section \ref{sec-SWGs}. Every time that the SWG from $U$ grows by a vertex, this vertex is paired to a half-edge incident to vertices in $\SWG^{\sss(1)}(\Toan+\vep)$ with positive probability, since there are of the order $n$ unpaired half-edges incident to $\SWG^{\sss(1)}(\Toan+\vep)$.

We claim that whp this path is not blocked by the competing type 2. Indeed, condition on $Q_n=\ell$. The vertices $\tilde U_1, \ldots, \tilde U_\ell$ found by the first-passage exploration from $\Ver$ are a size-biased reordering of the vertices in $[n]\setminus \SWG^{\sss(1,2)}(\Toan+\vep)$ (recall Definition \ref{def-size-biased-reor}). In particular, their degrees $d_{\tilde U_i}$ are tight random variables. 

Fix $i\in[\ell]$ and condition on $d_{\tilde U_i}=d$ for some $d\geq 1$. Then, since the passage times $(X^{\sss(2)}_e)_{e\in E(\CMnd)}$ are {\em independent} from $(X^{\sss(1)}_e)_{e}$, the probability that $\tilde U_i$ is occupied by type 2 is the same as that for any vertex $v$ with $d_v=d$. Since there are order $n$ of such vertices, and since $\SWG^{\sss(2)}(\Toan+b)=\op(n)$ for every $b>0$ by Proposition \ref{prop-exp-avoid-SWG} (recall also Corollary \ref{cor-losing-tight}), this implies that the probability that $\tilde U_i\in \SWG^{\sss(2)}(\Toan+b)$ vanishes. Since this is true for all $i\geq 1$ and all $d\geq 1$, in fact the probability that the first-passage percolation path from $U$ to $\SWG^{\sss(1)}(\Toan+\vep)$ is blocked by type 2 vanishes, as required.
\qed

\subsection{Winning type is hardly hindered by losing type}
\label{sec-winning-hardly-hindered-losing}
The next result describes how vertices are being found by the winning type $\mathcal{W}_n$ after time $T_{a_n/2}^{\sss(1)}\wedge T_{a_n/2}^{\sss(2)}=T_{a_n/2}^{\sss(\mathcal{W}_n)}$. The essence of our results is that $\Ntkwin$ and $\barLtwin$ develop in the same way as in a one-type process with type $\mathcal{W}_n$ {\em without competition}, as described in Theorem \ref{thm-degree-time-FPP}. 
We will see that at time $T_{a_n/2}^{\sss(\mathcal{W}_n)}+t$, a positive proportion of the vertices are found by the winning type. To describe how the winning type sweeps through the graph, we need some notation that adapts the notation used in Theorem \ref{thm-degree-time-FPP}. 

Write
	\eqn{
	\barNtkwin=\#\{v\colon D_v=k \mbox{ and } v \text{ is occupied by winning type at time } T_{a_n/2}^{\sss(\mathcal{W}_n)}+t\}/n.
	}
for the fraction of vertices that have degree $k$ and that have been occupied by the winning type at time $T_{a_n/2}^{\sss(\mathcal{W}_n)}+t$.
Further, for an edge $e=xy$ consisting of two half-edges $x$ and $y$ that are incident to vertices $v_x$ and $v_y$, we say that $e$ {\em spreads the winning type at time $s$} when $v_x$ (or $v_y$) is type $\mathcal{W}_n$ infected at time $s$, and $v_y$ (or $v_x$) is $\mathcal{W}_n$ infected at time $s$ through the edge $e$. We let
	\eqn{
	\Ltwin=\#\{e\colon \mbox{ $e$ has spread the winning infection at time }\Totan+t\}/[L_n/2],
	}
denote the fraction of edges that have spread the winning infection. Recall that $V^{\sss(i)}(k)$ denotes the conditional distribution of $V^{\sss(i)}$ conditionally on starting from a vertex of degree $k$, so that $V^{\sss(i)}$ has the same law as $V^{\sss(i)}(D).$ Further, conditionally on $\mathcal{W}=i$,  define $\Vwin(k)=V^{\sss(i)}(k).$ Recall that $D^\star$ denotes a size-biased version of a degree variable. Our result shows that, in the presence of competition, Theorem \ref{thm-degree-time-FPP} still holds for the winning type:

\begin{proposition}[Fraction of fixed degree winning type vertices and edges at fixed time]
\label{prop-degree-time}
As $n\rightarrow \infty$,
	\eqn{
	\label{Ntk-conv}
	\barNtkwin\convp \prob(\Vwin(k)\leq t)\prob(D=k),
	}
and
	\eqn{
	\label{Ltwin-conv}
	\barLtwin \convp \prob\Big(\Xwin+\tildeVwin_a\wedge\tildeVwin_b\leq t\Big),
	}
where, conditionally on $\mathcal{W}=i$, $(\tildeVwin_a,\tildeVwin_b)$ are two independent copies of $V^{\sss(i)}(D^\star-1)$, while $\Xwin=X^{\sss(i)}$.
\end{proposition}

We remark that Proposition \ref{prop-degree-time} can be interpreted by {\em first} conditioning on $\mathcal{W}_n=i$, and then considering the processes where $\mathcal{W}$ is everywhere replaced by $i$. Thus, while $\mathcal{W}_n$ is {\em random}, conditionally on $\mathcal{W}_n=i$, the random variables appearing in \eqref{Ntk-conv} and \eqref{Ltwin-conv} can be defined in terms of iid copies of $X^{\sss(i)}$, $D$ and $D^\star-1$.

The statement was proved in the first-passage percolation context in Theorem \ref{thm-degree-time-FPP}. We need to show that the presence of the losing type 2 infection started from vertex $\Ver_2$ does not affect this convergence result when $\mathcal{W}_n=1$. This part of the proof is an adaptation of the proof of Theorem \ref{thm-degree-time-FPP}, using ideas also appearing in the proof of Theorem \ref{thm-main}(a) in Section \ref{sec-proof-thm-a}.

\proof We start with \eqref{Ntk-conv}. We write
	\eqn{
	\barNtkwin=\barNtk-\bar{R}_n^{\sss(N)}(t),
	}
where $\barNtk$ is the proportion of vertices found by the first-passage percolation exploration of the winning type, where we ignore the competition with the losing type, and $\bar{R}_n^{\sss(N)}(t)$ corresponds to the proportion of vertices that are blocked by the losing type. We then use Theorem \ref{thm-degree-time-FPP} to obtain that 
	\eqn{
	\barNtk\convp \prob(\Vwin(k)\leq t)\prob(D=k).
	}
Thus, we are left to show that $\bar{R}_n^{\sss(N)}(t)\convp 0$. For this, we use Theorem \ref{thm-main}(a) and its proof, particularly Proposition \ref{prop-exp-avoid-SWG}. We use the first moment method, and write
	\eqn{
	\expec[\bar{R}_n^{\sss(N)}(t)]=\prob(U \text{ found by winning type at time $t$, and blocked by losing type}).
	}
Without loss of generality, we again assume that $\Toan<\Ttan$, so that type 1 is whp the winning type. Conditionally on this event, we know that $\SWG^{\sss(2)}(\Toan+\vep)=\SWG^{\sss(2)}(\Toan)$ occurs with probability close to one for $\vep>0$ small, and we assume this from now on. Denote 
	\eqn{
	\mathcal{E}_n^{\sss(1)}(\vep)=\{\Toan<\Ttan, \SWG^{\sss(2)}(\Toan+\vep)=\SWG^{\sss(2)}(\Toan)\},
	} 
so that, for $n\rightarrow \infty$ folllowed by $\vep \searrow 0$,
	\eqn{
	\expec[\bar{R}_n^{\sss(N)}(t)]=\expec[\indicwo{\mathcal{E}_n^{\sss(1)}(\vep)}\bar{R}_n^{\sss(N)}(t)]+o(1).
	}
Note that
	\eqan{
	&\{U \text{ found by winning type at time $\Ttan+t$, and blocked by losing type}\}\cap \mathcal{E}_n^{\sss(1)}(\vep)\nn\\
	&\quad 
	\subseteq \{W_{n}^{\sss (2)}(\Ver_2,U) \text{ in } [n]\setminus \SWG^{\sss(1)}(T_{a_n/2}^{\sss(1)}+\varepsilon)\leq T_{a_n/2}^{\sss(1)}+t\}
	\cap \mathcal{E}_n^{\sss(1)}(\vep),
	}
and, by Proposition \ref{prop-exp-avoid-SWG} together with the fact that $T_{a_n/2}^{\sss(1)}$ is tight, 
	\eqn{
	\expec[\indicwo{\mathcal{E}_n^{\sss(1)}(\vep)}\bar{R}_n^{\sss(N)}(t)]
	=\prob(\{W_{n}^{\sss (2)}(\Ver_2,U) \text{ in } [n]\setminus \SWG^{\sss(1)}(T_{a_n/2}^{\sss(1)}+\varepsilon)\leq T_{a_n/2}^{\sss(1)}+t\}
	\cap \mathcal{E}_n^{\sss(1)}(\vep))=o(1).
	}
This completes the proof of \eqref{Ntk-conv}.

The proof of \eqref{Ltwin-conv} is similar, where now we start from 
	\eqn{
	\barLtwin=\barLt -\bar{R}_n^{\sss(L)}(t),
	}
where $\barLtwin$ corresponds to the proportion of edges explored by the winning type first-passage percolation exploration where we ignore the competition with the losing type, and $\bar{R}_n^{\sss(L)}(t)$ corresponds to the proportion of such edges that are blocked by the losing type. 
Then $\barLt\convp \prob\Big(\Xwin+\tildeVwin_a\wedge\tildeVwin_b\leq t\Big)$ by Theorem \ref{thm-degree-time-FPP}, and we are left to show that 
$\bar{R}_n^{\sss(L)}(t)\convp 0.$ This is similar to the proof of \eqref{Ntk-conv}, recall the proof of Theorem \ref{thm-main}(a) in Section \ref{sec-proof-thm-a}, as well as Proposition \ref{prop-exp-avoid-SWG}.
\qed

\section{The winning type occupies almost all: Proof of Theorem \ref{thm-main}(b)}
\label{sec-almost-all}
In this section, we prove Theorem \ref{thm-main}(b). This proof is a minor modification of the proof of \cite[Theorem 1.1(b)]{DeiHof16} in \cite[Section 4]{DeiHof16}. Due to the fact that our estimates in Sections \ref{sec-FPP-prelim}--\ref{sec-most} are quite different from those in \cite[Section 3]{DeiHof16}, and the general edge-weight distributions do not allow for a simple formulation of the distribution of the explosion times $V^{\sss(i)}(k)$, we do present the entire proof here.

Throughout this section, we deal with the competition process, and explore the competition from the two vertices $\Ver_1$ and $\Ver_2$ simultaneously. As before, we write $\mathcal{W}_n$ for the type that occupies the largest number of vertices at time $T_{a_n/2}^{\sss(1)}\wedge T_{a_n/2}^{\sss(2)}=T_{a_n/2}^{\sss(\mathcal{W}_n)}$ and $\mathcal{L}_n$ for the type that occupies the smallest number of vertices. 

This section is organised as follows. In Section \ref{sec-losing-vertices-at-explosion}, we study the number of vertices occupied by the losing type at time $T_{a_n/2}^{\sss(\mathcal{W}_n)}$, and show that this converges in distribution jointly with $T_{a_n/2}^{\sss(\mathcal{W}_n)}$. In the remainder of the section, we prove that the asymptotic number $\Nlos^{**}$ of vertices occupied by type $\mathcal{L}_n$ after time $T_{a_n/2}^{\sss(\mathcal{W}_n)}$ also converges in distribution to an almost surely finite random variable. This is done in Sections \ref{sec-exploration-losing}--\ref{sec-conv-losing-vertices}, where Section \ref{sec-conv-losing-vertices} contains the proof of Theorem \ref{thm-main}(b). Section \ref{sec-exploration-losing} proves the weak convergence of the exploration process of the losing type, while Section \ref{sec-vertices-losing-large-times} studies the large-time asymptotics of the limiting process.

\subsection{Number of vertices occupied by losing type at explosion time}
\label{sec-losing-vertices-at-explosion}

Our first result is that $T_{a_n/2}^{\sss(\mathcal{W}_n)}$ converges to the minimum of the explosion times $V^{\sss(1)}$ and $V^{\sss(2)}$ of the one-type exploration processes, and that the asymptotic number $\Nlos^*$ of vertices that are then occupied by type $\mathcal{L}_n$ is finite. 

We start by introducing some notation. Recall that $\mathcal{W}$ and $\mathcal{L}$ denote the winning and the losing type, respectively, in the limit as $n\to\infty$. According to Theorem \ref{thm-main}(a), asymptotically type 1 wins with probability $\prob(V^{\sss(1)}<V^{\sss(2)})$ and type 2 with probability $\prob(V^{\sss(1)}>V^{\sss(2)})$. Then, $V^{\sss(1)}\wedge V^{\sss(2)}=\Vwin$ is close to the time when the winning type finds vertices of very high degree. The random variable $\Vlos$ does not have such a simple interpretation in terms of the competition process, since the winning type starts interfering with the exploration of the losing type before time $\Vlos$. The main aim of this section is to describe the exploration of the winning and losing types after time $\Vwin$, where the age-dependent branching process approximation breaks down and the types start interfering. The relation between the number of vertices found by the losing type and $\Vwin$ is described in the following lemma:
	
\begin{lemma}[Status at explosion of winning type]
\label{le:tt_explosion}
Let $\NLna(n)=\max\{m\colon T^{\sss (\sss\mathcal{L}_n)}_m \leq T_{a_n/2}^{\sss(\mathcal{W}_n)}\}$. Then, as $n\rightarrow \infty$,
	\eqn{
	(T_{a_n/2}^{\sss(\mathcal{W}_n)}, \NLna(n))\convd (\Vwin, \Nlos^*),
	}
where
	\begin{equation}
	\label{Nlos-def}
	\Nlos^*\stackrel{d}{=}\max\big\{m\colon \Tlos_m \leq \Vwin\big\}.
	\end{equation}
Here, $(\Tlos_m)_{m\geq 0}$ are the vertex-discovery times in the unimodular $(D,X^{\sss(\mathcal{L})})$ branching process conditioned on $\Vwin<\Vlos$ (recall Section \ref{sec-SWGs}). Since $\Vwin<\Vlos$ holds a.s., the random variable $\Nlos^*$ is finite a.s.
\end{lemma}

\proof
By definition, the number of vertices occupied by type $\mathcal{W}_n$ at time $T_{a_n/2}^{\sss(\mathcal{W}_n)}$ equals $a_n/2$. Furthermore, by Proposition \ref{prop-prel}(c), the set of type 1 and type 2 occupied vertices, respectively, are whp disjoint at this time, that is, none of the infection types has then tried to occupy a vertex that was already occupied by the other type. Up to that time, the exploration processes started from vertex 1 and 2, respectively, hence behave like in the corresponding one-type processes. The asymptotic distributions of $T_{a_n/2}^{\sss(\mathcal{W}_n)}$ and $\NLna(n)$ follow from the characterisation of the times $T_m^{\sss(1)}$ and $T_m^{\sss(2)}$ in the two one-type processes described in Section \ref{sec-SWGs}, and the convergence result in Proposition \ref{prop-prel}(c).
\qed

\subsection{Convergence of the exploration process of losing type}
\label{sec-exploration-losing}
We grow the SWG of type $\mathcal{L}_n$ from size $\NLna(n)$ onwards. At this moment, \whp the type $\mathcal{L}_n$ has not yet tried to occupy a vertex that was already occupied by type $\mathcal{W}_n$. However, when we grow the SWG further, then type $\mathcal{W}_n$ will grow very quickly due to its explosion. We will show that the growth of type $\mathcal{L}_n$ is thus delayed to the extent that it will only occupy {\em finitely} many vertices. An important tool in proving this rigorously is a stochastic process $(S_m')_{m\geq 0}$ keeping track of the number of unexplored half-edges incident to the SWG of the losing type.

Recall that, by the construction of the two-type exploration process described in Sections \ref{sec-SWGs} and \ref{sec-prel}, the quantity $S_0'(n)$ represents the number of (unpaired) half-edges incident to the SWG of type $\mathcal{L}_n$ at time $T_{a_n/2}^{\sss(\mathcal{W}_n)}$. By Lemma \ref{le:tt_explosion}, 
	\eqn{
	\label{S_0'-def}
	S_0'(n)\convd S_0'\equiv D^{\sss (\mathcal{L})}+\sum_{i=2}^{\Nlos^{*}} B_i^{\sss (\mathcal{L})}-\Nlos^{*},
	}
where $D^{\sss (\mathcal{L})}$ has distribution $D$, and $(B_i^{\sss (\mathcal{L})})_{i\geq 1}$ are iid copies of $D^\star-1$. Conditionally on $\mathcal{W}=i$, denote the remaining life-times of the unimodular $(D,X^{\sss(i)})$ process by $(R_{0,j}^{\sss(\mathcal{L}_n)})_{j\leq S_0'(n)}$. By Lemma \ref{le:tt_explosion} and its proof, and conditionally on $\Vwin$,
	\eqn{
	\label{R_0'-def}
	(R_{0,j}^{\sss(\mathcal{L}_n)})_{j\leq S_0'(n)}\convd (R_{0,j}')_{j\leq S_0'},
	}
where, conditionally on $\mathcal{W}=i$, $S_0'=\ell$ and $\Vwin=x$,  $(R_{0,j}')_{j\leq S_0'}$ are the conditional remaining lifetimes of the unimodular $(D,X^{\sss(3-i)})$ process at time $x$ and conditionally on $S_0'=\ell$ and $V^{\sss(3-i)}>x$. We also define $T_m'(n)=0$.

The sequences $(T_m'(n))_{m\geq 0}, (S_m'(n))_{m\geq 0}$ and $(R_{m,j}^{\sss(\mathcal{L}_n)})_{j\leq S_m'(n)}$ are constructed recursively and simultaneuously.  We grow the SWG of type $\mathcal{L}_n$ one edge at a time by pairing the half-edge with minimal remaining edge-weight to a uniform unpaired half-edge. Denote the half-edge of minimal weight in the $m$th step by $x_m$ and the half-edge to which it is paired by $P_{x_m}$, where we recall that $v_y$ denotes the vertex incident to the half-edge $y$. Let the remaining edge-weight of the half-edge $x_m$ be denoted by $\underline{R}_m(n)$. We let $T_{m+1}'(n)=T_m'(n)+\underline{R}_m(n)$. Of course, it is possible that  $v_{P_{x_m}}$ is already infected by the winning or losing type, and then the SWG of the losing type does not grow and we define $S_{m+1}'(n)=S_m'(n)-1$ and $(R_{m+1,j}^{\sss(\mathcal{L}_n)})_{j\leq S_{m+1}'(n)}=(R_{m,j}^{\sss(\mathcal{L}_n)}-\underline{R}_m(n))_{j\leq S_m'(n)}$ (where the one value that is equal to zero is removed from the list). When  $v_{P_{x_m}}$ is not yet infected by the winning or losing type, then we let $B_m'(n)$ denote $D_{U_{P_{x_m}}}-1$. Let
	\begin{equation}
	\label{eq:Sm_rec}
	S_{m+1}'(n)-S_{m}'(n)=B_{m+1}'(n)-1.
	\end{equation}
When $B_{m+1}'(n)\geq 1$, we draw $B_{m+1}'(n)$ iid random variables with law $X^{\sss(i)}$, conditionally on $\mathcal{W}=i$ and call these $\big(R_{m+1,j}^{\sss(\mathcal{L}_n)}\big)_{j=1}^{B_m'(n)}$. We then let $\big(R_{m+1,j}^{\sss(\mathcal{L}_n)}\big)_{j=1}^{S_{m+1}'(n)}$ denote the combined list $\big(R_{m+1,j}^{\sss(\mathcal{L}_n)}\big)_{j=1}^{B_m'(n)}$ together with $(R_{m,j}^{\sss(\mathcal{L}_n)}-\underline{R}_m(n))_{j\leq S_m'(n)}$, where the unique component equal to zero is removed, so that all coordinates are strictly positive.

Our aim is to identify the scaling limit of $(T_m'(n),S_m'(n), (R_{m,j}^{\sss(\mathcal{L}_n)})_{j\leq S_m'(n)})_{m\geq 0}$. To this end, we define $S_0'$ as in \eqref{S_0'-def}, $T_0'=0$ and 
$(R_{0,j}')_{j\leq S_0'}$ as in the rhs of \eqref{R_0'-def}. Further, for $m\geq 1$, again define $(T_m',S_m')_{m\geq 0}$ recursively by $T_{m+1}'-T_{m}'=\underline{R}_m=\min_{i\leq S_m'} R_{m,j}'$. Further, $(S_m')_{m\geq 0}$ is defined recursively by $S_{m+1}'-S_{m}'=B_{m+1}'-1$, where, conditionally on $T_{m}'=t$, for all $k\geq 1$,
	\eqn{
	\label{Bm'-def}
	\prob(B_{m+1}'=k\mid T_{m+1}'=t)
	=\prob(D^\star=k+1)\frac{\prob(\Vwin(k+1)>t)}
	{\prob\Big(\Xwin+\tildeVwin_a\wedge\tildeVwin_b> t\Big)},
	}
while $\prob(B_{m+1}'=0\mid T_{m+1}'=t)=1-\sum_{k\geq 1} \prob(B_{m+1}'=k\mid T_{m+1}'=t)$ corresponds to the asymptotic probability that the vertex that the losing type is trying to occupy has already been occupied by the winning type. Further, $(R_{m+1,j}^{\sss(\mathcal{L})})_{j\leq S_m'}$ is a combination of the lists $\big(R_{m+1,j}^{\sss(\mathcal{L}_n)}\big)_{j=1}^{B_m'(n)}$ and $(R_{m,j}^{\sss(\mathcal{L})}-\underline{R}_m)_{j\leq S_m'}$, where the unique component equal to zero is removed, so that all coordinates are strictly positive. We denote this process by $(T_m',S_m', (R_{m,j}^{\sss(\mathcal{L})})_{j\leq S_m'})_{m\geq 0}$.

\begin{remark}[Edge-weight distribution vs.\ weights on half-edges]
\label{rem-edge-weight} 
{\rm In the above construction, we explore from vertices 1 and 2 simultaneously, and search for the minimal weight among unexplored half-edges of the loosing type. This half-edge is paired to a randomly chosen second half-edge, and the passage time of the resulting edge should be given by the weight of this half-edge, that is, the second half-edge should not be assigned any weight at all. The careful reader may note that this is not the case in the above construction when the second half-edge belongs to a vertex that is already infected. In that case, however, the edge will never be used to transmit infection (indeed, such edges are not included in the SWG defined in Section \ref{sec-SWGs}) and its assigned passage time is hence unimportant for the competition process.}\hfill $\ensymboldefinition$
\end{remark}
	
The following lemma shows that $(T_l',S_l', (R_{l,j}^{\sss(\mathcal{L})})_{j\leq S_l'})_{l\geq 0}$ is indeed the limit in distribution of the process $(T_l'(n),S_l'(n), (R_{l,j}^{\sss(\mathcal{L}_n)}(n))_{j\leq S_l'(n)})_{l\geq 0}$. In its statement, we recall that $\field_{a_n}$ denotes the $\sigma$-field of the exploration of the two competing species up to time $T_{a_n/2}^{\sss(\mathcal{W}_n)}$:

\begin{lemma}[Exploration of the losing type beyond explosion of the winning type]
\label{lem-expl-losing}
Conditionally on $\field_{a_n}$, and for all $m\geq 1$, as $n\rightarrow \infty$,
	\eqn{
	(T_l'(n),S_l'(n), (R_{l,j}^{\sss(\mathcal{L}_n)}(n))_{j\leq S_l'(n)})_{l=0}^m\convd (T_l',S_l', (R_{l,j}^{\sss(\mathcal{L})})_{j\leq S_l'})_{l=0}^m.
	}
\end{lemma}

\proof We prove the claim by induction on $m$. The statement for $m=0$ follows from $T_0'(n)=T_0'=0$, Lemma \ref{le:tt_explosion} and \eqref{S_0'-def}--\eqref{R_0'-def}. This initialises the indiction hypothesis.

To advance the induction hypothesis, we introduce some further notation. Let $\field_m'$ be the $\sigma$-field generated by $\field_{a_n}$ together with $(T_l'(n),S_l'(n), (R_{l,j}^{\sss(\mathcal{L}_n)}(n))_{j\leq S_l'(n)})_{l=0}^m$. Then, conditionally on $\field_m'$,
	\eqn{
	T_{m+1}'(n)-T_{m}'(n)=\underline{R}_m(n)=\min_{j\leq S_m'(n)} R_{m,j}^{\sss(\mathcal{L}_n)}(n)\convd \min_{j\leq S_m'} R_{m,j}^{\sss(\mathcal{L})}=\underline{R}_m,
	}
by the fact that $(T_l'(n),S_l'(n), (R_{l,j}^{\sss(\mathcal{L}_n)}(n))_{j\leq S_l'(n)})_{l=0}^m\convd (T_l',S_l, (R_{l,j}^{\sss(\mathcal{L})})_{j\leq S_l'(n)})_{l=0}^m$ by the induction hypothesis. This advances the claim for $T_{m}'(n)$. For $S_{m+1}'(n)$, we note that $B_{m+1}'(n)=k$ precisely when the half-edge that is found is paired to a half-edge of a vertex of degree $k+1$ that is not yet infected. The number of vertices that is type $\mathcal{L}_n$-infected is negligible. Therefore, writing $N_n^{\sss(k+1)}$ for the total number of vertices of degree $k+1$, and noting that at time $t$, there are 
$L_n-2\Ltwin+\op(n)$ many unpaired half-edges,
	\eqn{
	\prob(B_{m+1}'(n)=k\mid \field'_{T_{m+1}'(n)}, T_{m+1}'(n)=t)
	=\frac{(k+1)[N_n^{\sss(k+1)}-\Ntkpluswin]}{L_n-2\Ltwin}(1+\op(1)).
	}
We rewrite this as
	\eqn{
	\prob(B_{m+1}'(n)=k\mid \field'_{T_{m}'(n)}, T_{m+1}'(n)=t)=\frac{(k+1)}{(L_n/n)}\frac{\bar{N}_n^{\sss(k+1)}-\barNtkpluswin}{1-2\barLtwin}(1+\op(1)),
	}
where $\bar{N}_n^{\sss(k+1)}=N_n^{\sss(k+1)}/n$ denotes the proportion of vertices with degree $k+1$.
By Proposition \ref{prop-degree-time}, this is equal to
	\eqan{
	&\prob(B_{m+1}'(n)=k\mid \field'_{T_{m}'(n)}, T_{m+1}'(n)=t)\\
	&\qquad \convp\frac{(k+1)}{\expec[D]}\frac{\prob(D=k+1)-\prob(\Vwin(k+1)\leq t)\prob(D=k+1)}{1-\prob\Big(\Xwin+\tildeVwin_a\wedge\tildeVwin_b\leq t\Big)}\nn\\
	&\qquad= \prob(D^\star=k+1)\frac{\prob(\Vwin(k+1)>t)}
	{\prob\Big(\Xwin+\tildeVwin_a\wedge\tildeVwin_b\big)> t\Big)}=\prob(B_{m+1}'=k\mid T_{m+1}'=t),\nn
	}
as required. This shows that, conditionally on $\field'_{T_{m}'(n)}$ and $T_{m+1}'(n)=t$, the law of $(T_{m+1}'(n)-T_{m}'(n), B_{m+1}'(n))$ converges to \eqref{Bm'-def}. In both constructions, we add $B_{m+1}'(n)$ and $B_{m+1}'$ iid random edge-weights to the lists $(R_{m,j}^{\sss(\mathcal{L})}(n))_{j\leq S_m'(n)}$ and $(R_{m,j}^{\sss(\mathcal{L})})_{j\leq S_m'})$, subtract $\underline{R}_m(n)$, respectively $\underline{R}_m$, from the existing elements in the list, and remove the unique zero. As a result, also $(R_{m+1,j}^{\sss(\mathcal{L})}(n))_{j\leq S_{m+1}'(n)}\convd (R_{l,j}^{\sss(\mathcal{L})})_{j\leq S_m'})$, conditionally on $\field'_{T_{m}'(n)}, T_{m+1}'(n)-T_{m}'(n)$ and $B_{m+1}'(n)$. This advances the induction hypothesis, and thus completes the proof.\qed\medskip

Denote $H'(n)=\min\{m\colon S_m'(n)=0\}$ and $H'=\min\{m\colon S_m'=0\}$. In the following corollary, we show that $H'(n)\convd H'$.

\begin{corollary}[Convergence of hitting of zero]
\label{lem-hitting-time-zero}
For all $m\geq 1$, as $n\rightarrow \infty$,
	\eqn{
	\prob(H'(n)\leq m\mid \field_{a_n})\convp \prob(H'\leq m\mid \field_{\Vwin}).
	}
Therefore, $H'(n)\convd H'$, where $H'$ is possibly defected (i.e., $\prob(H'=\infty)\in [0,1]$).
\end{corollary}

\proof It suffices to realise that the event $\{H'(n)\leq m\}$ is measurable wrt the sequence $(T_l'(n),S_l'(n), (R_{l,j}^{\sss(\mathcal{L}_n)}(n))_{j\leq S_l'(n)})_{l=0}^m$. Then the claim follows from Lemma \ref{lem-expl-losing}.
\qed

\subsection{Asymptotics of vertices occupied by losing type for large times}
\label{sec-vertices-losing-large-times}
Note that $B_m'$ in (\ref{Bm'-def}) has infinite mean when we condition on $T_m'=t=0$, which implies that initially many of its values are large. This is the problem that we need to overcome in showing that the number of vertices found by the losing type is finite. As it turns out, conditionally on $T_m'=t$, the mean of $B_m'$ decreases as $t$ increases, and becomes smaller than 1 for large $t$, so this saves our day. In order to prove this, we first need some results on the process $\Vwin(k)$:

\begin{lemma}[Bounds and asymptotics for $\Vwin(k)$]\label{lem-V(k)-asymp}\leavevmode
The law of $\Vwin(k)$ is related to that of $\Vwin(1)$ by
	\eqn{
	\label{V(k)-rec}
	\prob(\Vwin(k)>t)=\prob(\Vwin(1)>t)^k,\qquad k\geq 1, t\geq 0.
	}
Further, for all $t\geq 0$,
	\eqn{
	\label{EVaVb_bd}
	\prob\left(\Xwin+\tildeVwin_a\wedge\tildeVwin_b> t\right)\geq \prob(\Vwin(1)>t)^2.
	}
\end{lemma}
	
\proof The relation \eqref{V(k)-rec} follows since $\Vwin(k)$ is the explosion time starting from $k$ individuals, which is the minimum of the explosion times of $k$ iid explosion times starting from 1 individual, i.e.,
	$$
	\Vwin(k)\stackrel{d}{=}\min_{i=1}^k \Vwin_i,
	$$
where $(\Vwin_i)_{i\geq 1}$ are iid with law $\Vwin(1)$ and $\Xwin$ is the edge-weight of the winning type. From this, \eqref{V(k)-rec} follows immediately. For \eqref{EVaVb_bd}, write $G(t)=\prob(\Vwin(1)>t)$ and note that $\Vwin(1)\stackrel{d}{=}\Xwin+\min_{i=1}^{D^{\star}-1} \Vwin_i$, where again $(\Vwin_i)_{i\geq 1}$ are iid with law $\Vwin(1)$. Thus, conditioning on $\Xwin$ and $D^{\star}$, and using \eqref{V(k)-rec}, leads to
	\eqn{
	\label{V(1)-rec}
	G(t)=[1-F_{\sss X}](t)+\int_0^t f_{\sss X}(s) \expec[G(t-s)^{D^{\star}-1}]ds,
	}
where $[1-F_{\sss X}](t)=\prob(\Xwin>t)$. Furthermore, since $\tildeVwin_a$ and $\tildeVwin_b$ are iid with the same distribution as $\Vwin(D^{\star}-1)\stackrel{d}{=}\min_{i=1}^{D^{\star}-1} \Vwin_i$, we similarly obtain that
	\eqn{
	\label{edge_rec}
	\prob\left(\Xwin+\tildeVwin_a\wedge\tildeVwin_b> t\right)=[1-F_{\sss X}](t)+ \int_0^t f_{\sss X}(s)\expec[G(t-s)^{D^{\star}-1}]^2ds.
	}
By the Cauchy-Schwarz inequality (where we split $f_{\sss X}(s)=\sqrt{f_{\sss X}(s)}\sqrt{f_{\sss X}(s)}$ in the lhs),
	\eqn{
	\label{Cauchy}
	\left(\int_0^t f_{\sss X}(s)\expec[G(t-s)^{D^{\star}-1}]ds\right)^2\leq F_{\sss X}(t)\int_0^t f_{\sss X}(s) \expec[G(t-s)^{D^{\star}-1}]^2ds.
	}
Combining \eqref{V(1)-rec}, \eqref{edge_rec} and \eqref{Cauchy} yields that
	\eqan{
	\prob\left(\Xwin+\tildeVwin_a\wedge\tildeVwin_b> t\right)&\geq [1-F_{\sss X}](t)+\frac{(G(t)-[1-F_{\sss X}](t))^2}{F_{\sss X}(t)}.
	}
Working out the square gives
	\eqan{
	\prob\left(\Xwin+\tildeVwin_a\wedge\tildeVwin_b> t\right)&\geq [1-F_{\sss X}](t)+\frac{[1-F_{\sss X}](t)^2}{F_{\sss X}(t)}
	-2\frac{G(t)[1-F_{\sss X}](t)}{F_{\sss X}(t)}+\frac{G(t)^2}{F_{\sss X}(t)}.
	}
We use that
	\eqn{
	[1-F_{\sss X}](t)+\frac{[1-F_{\sss X}](t)^2}{F_{\sss X}(t)}=[1-F_{\sss X}](t)\Big[1+\frac{[1-F_{\sss X}](t)}{F_{\sss X}(t)}\Big]=\frac{[1-F_{\sss X}](t)}{F_{\sss X}(t)},
	}
and
	\eqn{
	\frac{G(t)^2}{F_{\sss X}(t)}=G(t)^2 +\frac{G(t)^2[1-F_{\sss X}(t)]}{F_{\sss X}(t)}.
	}
Rearranging terms then gives
	\eqan{
	\prob\left(\Xwin+\tildeVwin_a\wedge\tildeVwin_b> t\right)&\geq G(t)^2+\frac{[1-F_{\sss X}(t)](G(t)-1)^2}{F_{\sss X}(t)}\geq \, G(t)^2,
	}
as desired.
\qed
\medskip

Recall that $T_m'$ denotes the minimal residual lifetime of a half-edge in the losing exploration, after pairing $m$ half-edges. Let $B_m'(t)\stackrel{d}{=}B'_m|\,T'_m=t$. The next lemma shows that $(B'_m(t))_{m\geq 1}$ is stochastically dominated by an iid sequence whose mean is strictly smaller than 1 for large $t$. It also shows that $T'_m\to\infty$ almost surely. It is a key ingredient in the proof of Theorem \ref{thm-main}(b):

\begin{lemma}[Asymptotic behavior of $B_m'(t)$ and $T'_m$]
\label{le:B'm_dom}
\begin{itemize}
\item[\rm{(a)}] For each fixed $t>0$, the sequence $(B'_m(t))_{m\geq 1}$ is stochastically dominated by an iid sequence $(\bar{B}(t))_{m\geq 1}$. Furthermore, $\expec[\bar{B}(t)]$ is finite for each fixed $t>0$ and $\expec[\bar{B}(t)]\to \prob(D^{\star}=2)$ as $t\to\infty$.
\item[\rm{(b)}] $T'_m\to\infty$ almost surely as $m\to\infty$.
\end{itemize}
\end{lemma}

\proof Recalling the definition of $B'_m|\,T'_m=t$ in \eqref{Bm'-def}, using both \eqref{V(k)-rec} and \eqref{EVaVb_bd} in Lemma \ref{lem-V(k)-asymp}(a), we obtain for $k\geq 1$ that
	\eqan{
	\prob(B'_m(t)=k) & =\prob(D^\star=k+1)\frac{\prob(\Vwin(k+1)>t)}
	{\prob\Big(\Xwin+\tildeVwin_a\wedge\tildeVwin_b> t\Big)}\nn\\
	&\leq \prob(D^\star=k+1)\prob(\Vwin(1)>t)^{k-1}.
	}
	
Denote $\bar{p}_k(t)=\prob(D^\star=k+1)\prob(\Vwin(1)>t)^{k-1}$, and let $\bar{B}(t)$ be defined by
	\eqn{
	\prob(\bar{B}(t)=k)=\left\{
	\begin{array}{ll}
	\bar{p}_k(t) & \mbox{for }k\geq 1;\\
	1-\sum_{k\geq 1} \bar{p}_k(t) & \mbox{for }k=0.
  	\end{array}
            \right.
	}
Then, obviously, $B_m(t)$ is stochastically dominated by $\bar{B}(t)$.
For any fixed $t>0$, we have that $\prob(\Vwin(1)>t)<1$, so that $\bar{p}_k(t)\to 0$ exponentially in $k$. It follows that $\bar{B}(t)$ has all moments so that, in particular, its mean is finite. Furthermore, $\bar{p}_1(t)=\prob(D^\star=2)$ for any $t>0$, while $\bar{p}_k(t)\to 0$ as $t\to\infty$ for each $k\geq 2$. Hence $\expec[\bar{B}(t)]\rightarrow \prob(D^\star=2)$.

As for (b), recall the construction of the processes $(S'_m)_{m\geq 0}$ and $(T'_m)_{m\geq 0}$, with $S'_m=S'_0+\sum_{i=1}^mB'_i$, where $B_m'=B_m(T_{m}')$. Note that $B'_i(t)$ is stochastically bounded by an iid sequence that is decreasing in $t$ having finite mean for all $t>0$, and that $t=T_m'\geq T_1'>0$ a.s. This implies that, conditionally on $T_1'=t_1'>0$, the (inhomogeneous) exploration tree constructed by $(T_m',S_m', (R_{m,j}^{\sss(\mathcal{L})})_{j\leq S_m'})_{m\geq 0}$ can be bounded from above by an $(S_1', \bar{B}(t_1), X)$ age-dependent branching process with offspring distribution $\bar{B}(t_1)$, except for the root, which has offspring $S_1'$. This age-dependent branching process is conservative since $\expec[ \bar{B}(t_1)]<\infty$.  Thus, $T'_m\to\infty$ a.s.
\qed

\subsection{Convergence of the number of vertices occupied by the losing type: Proof of Theorem \ref{thm-main}(b)}
\label{sec-conv-losing-vertices}
Recall the construction of the process $(S_m')_{m\geq 0}$ in the recursion (\ref{eq:Sm_rec}), and recall that $\NLna(n)$ denotes the number of vertices infected by the losing type at time $T_{a_n/2}^{\sss(\mathcal{W}_n)}$. Denote the total number of vertices infected by type $\mathcal{L}_n$ {\em after} time $T_{a_n/2}^{\sss(\mathcal{W}_n)}$ by $\NLnaa(n)$. We can identify this as
	\eqn{
	\NLnaa(n)=\#\{m\colon B_m'(n)\geq 1\}.
	}
Indeed, each time when a new vertex is found that is not infected by type $\mathcal{W}_n$, the degree of the vertex is at least 2, so that $B_m'(n)\geq 1$. Thus, the number of vertices found is equal to the number of $m$ for which $B_m'(n)\geq 1$.

Recall that the total asymptotic number of losing type vertices is denoted by $\Nlos$. This number can now be expressed as
    	\eqn{
	\label{Nlos-tot-def}
	\Nlos=\Nlos^*+\Nlaa,
	}
where $\Nlos^*$ is defined in Lemma \ref{le:tt_explosion} and $\NLnaa(n)\convd \Nlaa :=\#\{m\colon B_m'\geq 1\}$, where the weak convergence follows from Lemma \ref{lem-expl-losing}. Further, since the convergence in Lemma \ref{lem-expl-losing} is {\em conditional on} $\field_{a_n},$ we also obtain the joint convergence
	\eqn{
	(\NLna(n),\NLnaa(n))\convd (\Nlos^*,\Nlaa),
	}
which implies \eqref{Nlos-def}. To prove Theorem \ref{thm-main}(b), it hence suffices to show that the random variable $\Nlaa$ is finite almost surely. This certainly follows when $H'=\max\{m\colon S_m'\geq 1\}$ is almost surely finite, which is what we prove below.

We argue by contradiction. Assume that $H'=\infty$. Then, $T_{m-1}'<\infty$ for every $m$.
Furthermore, $S'_m-S'_{m-1}=B'_m-1$ where, by Lemma \ref{le:B'm_dom}(a), the contribution $B'_m$ is stochastically dominated by $\bar{B}(T'_m)$, with $\expec[\bar{B}(t)]\to \prob(D^\star=2)$ as $t\to\infty$. Pick $k$ large so that $\expec[\bar{B}(T'_k)\mid T'_k]<1$, which is possible since $\prob(D^\star=2)<1$ and $T'_k\to\infty$ a.s.\ by Lemma \ref{le:B'm_dom}(b). Then, conditionally on $T'_k$ and for $m>k$, $S'_m-S'_k$ is stochastically dominated by $\sum_{i=k+1}^m(\bar{B}_i(T'_k)-1)$ -- a sum of iid variables with negative mean. It follows that $S'_m$ hits 0 in finite time, so that $H'<\infty$, which is a contradiction.
\qed

\paragraph{\bf Acknowledgement} The author thanks Enrico Baroni, Mia Deijfen and J\'ulia Komj\'athy for helpful discussions at the start of the project, and J\'ulia Komj\'athy for comments on an early draft. This work is supported in parts by the Netherlands Organisation for Scientific Research (NWO) through VICI grant 639.033.806 and the Gravitation {\sc Networks} grant 024.002.003

\bibliographystyle{abbrv}
\bibliography{../../bib/bib}

\end{document}